%% file: COMET.tex
\documentclass[journal,12pt,draftclsnofoot,onecolumn]{IEEEtran}
\input{input.tex}

\input{header.tex}

\usepackage{graphicx,epstopdf}
\usepackage{caption,subcaption}

\usepackage[ruled,vlined]{algorithm2e}
\usepackage{algorithmic}
\usepackage{xr-hyper}
\usepackage{hyperref}

\UseRawInputEncoding

\title{A New Class of Composite Objective Multi-step Estimating-sequence Techniques (COMET)}
\author{{Endrit~Dosti, Sergiy~A.~Vorobyov, and~Themistoklis~Charalambous}
\thanks{E. Dosti and S. A. Vorobyov are with the Department of Signal Processing and Acoustics, Aalto University, Finland, e-mails:  firstname.lastname@aalto.fi.}
\thanks{T. Charalambous is with the Department of Electrical Engineering and Automation, Aalto University, Finland, e-mail: firstname.lastname@aalto.fi.}}

\begin{document}
\maketitle

\begin{abstract}
We devise a new accelerated gradient-based estimating sequence technique for solving large-scale optimization problems with composite structure. More specifically, we introduce a new class of estimating functions, which are obtained by utilizing a tight lower bound on the objective function. Then, by exploiting the coupling between the proposed estimating functions and the gradient mapping technique, we construct a class of composite objective multi-step estimating-sequence techniques (COMET). We propose an efficient line search strategy for COMET, and prove that it enjoys an accelerated convergence rate. The established convergence results allow for step size adaptation. Our theoretical findings are supported by extensive computational experiments on various problem types and datasets. Moreover, our numerical results show evidence of the robustness of the proposed method to the imperfect knowledge of the smoothness and strong convexity parameters of the objective function.
\end{abstract}

\IEEEpeerreviewmaketitle

\vspace{-3mm}

\section{Introduction}
\label{sec:intro}
In this work, we devise accelerated black-box methods for solving large-scale convex optimization problems with a composite objective structure by using only first-order information. Such problems typically have the following structure:
\begin{equation}
	\begin{aligned}
		\label {opt_prob}
		& \underset{x \in \mathcal{R}^n}{\text{minimize}}
		& & F(x) = f(x) + \tau g(x), \tau > 0,
	\end{aligned}
\end{equation}
where the function $f : \mathcal{R}^n \rightarrow \mathcal{R}$ is an  $L_f$-smooth and $\mu_f$-strongly convex function with parameters such that $0 \leq \mu_f \leq L_f$. The regularizer $g: \mathcal{R}^n \rightarrow \mathcal{R}$ is a simple convex lower semi-continuous function with strong convexity parameter $\mu_g$. Typically, in signal processing applications the function $g(x)$ is simple, meaning that we are able to find a closed-form solution for minimizing the summation of $g$ and some simple auxiliary functions \cite{Nesterov_2007}. In more practical terms, the assumption on the simplicity of $g$ implies that its proximal map, defined as
\begin{equation}
	\begin{aligned}
		\label {prox_g}
		\text{prox}_{\tau g} \triangleq \; & \text{arg} \underset{z \in {\mathcal{R}^n}}{\min}
		& &\left(g(z) + \frac{1}{2 \tau} ||z - x||^2\right), x \in \mathcal{R}^n, 
	\end{aligned}
\end{equation}
is computed with complexity $\mathcal{O} (n)$ \cite{Boyd_prox_alg}. Here $|| \cdot||$ denotes the standard $l_2$ norm.

Problems that share the same structure as \eqref{opt_prob} arise quite often in different scientific disciplines, such as signal and image processing, data analysis and machine learning. Typical applications in which the formulation given in \eqref{opt_prob} is relevant include compressive sensing, phase retrieval problems, medical imaging, dictionary learning and many more (see \cite{Cevher_sp_mag, 6879577, 8119874, 9028239, gu2017projected} and references therein). When considering applications, the variable $x$ represents the model parameters, whereas the role of $f (x)$ is to ensure a good fit between the observed data and the estimated parameters. In signal processing applications, $g (x)$ acts as a regularizer and typically takes the form of some parameter shrinkage norm, i.e., $l_2$ norm \cite{9018066, Wainwright}, sparsity-enforcing norm, i.e., $l_1$ norm \cite{tibshirani1996regression, tropp2010computational, combettes2011proximal}, or its counterpart for the rank function, i.e., the nuclear norm \cite{candes2015phase, yurtsever2015scalable}. Another popular structure for $g(x)$ is the Chebyshev norm, i.e., the $l_\infty$ norm \cite{demo}. The function $g(x)$ can also be used to embed convex constraints, in which case it would act as an indicator function of some closed convex set \cite{Nesterov_2007}. 

In the context of large-scale optimization \cite{nesterov2014subgradient}, problems that share the same structure as \eqref{opt_prob} are solved iteratively using different first-order optimization algorithms \cite{Beck_book, d2021acceleration}. The bounds on the performance of black-box first-order methods have been established by Nemirovskiy and Yudin in \cite{Nemirovski_Yudin}. Loosely speaking, a first-order method is optimal in the black-box framework if it achieves the accelerated convergence rate with respect to the iteration counter $k$, while at the same time complying with the lower complexity bounds. The question of how to construct practical methods that are optimal has attracted the attention of the research community over the decades. One of the first methods that managed to achieve the accelerated convergence rate in the black-box framework was the heavy ball method \cite{Polyak}. Therein, the acceleration is achieved by adding a momentum term to the gradient step, which nudges the new iterate in the direction of the previous step. The first method that is optimal in the sense of \cite{Nemirovski_Yudin} is the fast gradient method (FGM) \cite{Nesterov_83}. The method is built by using the mathematical machinery of estimating sequences, and has been since widely studied in the literature (see for instance \cite{Auslender_Teboulle, Lan2011PrimaldualFM, candes_donoghue, Aspremont, schmidt2011convergence, devolder2014first}).

A large portion of the recent research in first-order optimization is also focused in finding different reasons behind acceleration. For instance, in \cite{Allen-Zhu} the authors have constructed accelerated first-order methods by exploiting the linear coupling between mirror and gradient descent. The authors of \cite{Bubeck_geod} have derived an accelerated first-order method, which was inspired by the ellipsoid method. The proposed method is efficient, however, it suffers from the drawback that it requires an exact line search. An interesting framework is the one established in \cite{pmlr-v40-Flammarion15, Su_Boyd_Candes, WibisonoE7351}, wherein the authors model the continuous-time limit of FGM as a second-order differential equation. Then, FGM equations can be obtained based on such a framework. Moreover, the convergence rate for FGM can also be derived by using theory from robust control \cite{Recht}. A novel approach for analyzing the worst-case performance of first-order black-box optimization methods has appeared in \cite{drori-teboulle}. The analysis conducted therein relies on the observation that the worst-case behavior improvement of a black-box method is itself an optimization problem, which is referred to as the Performance Estimation Problem. By utilizing this approach, the authors of \cite{kim_fessler_1} and \cite{kim_fessler_2} introduce optimized first-order methods that are computationally efficient and achieve a convergence bound that is two times smaller than the one attained by FGM. However, the development of these algorithms is only restricted to solving optimization problems with differentiable objective functions.  

Among the various approaches to the acceleration of first-order methods that were discussed above, the methods that were built based on the machinery of estimating sequences have attracted a lot of attention (see \cite{Bubeck_book, d2021acceleration} and references therein). Let us now highlight several reasons that have led to their success. First, on a theoretical level FGM-type methods are proven to be optimal in the sense of \cite{Nemirovski_Yudin}. Second, their practical performance is competitive even when they are used in conjunction with simple line search strategies such as backtracking\cite{Iulian_1, nesterov2015universal}. Third, they can be scaled to construct accelerated second order methods \cite{fast_newton_method} and accelerated higher order methods \cite{nesterov2020inexact}. Lastly, they have been shown to excel in performance even when they have been extended to other settings such as distributed optimization \cite{jakovetic2014fast}, nonconvex optimization \cite{ghadimi2016accelerated}, stochastic optimization \cite{kulunchakov:hal-01993531}, non-Euclidean optimization \cite{ahn2020nesterov, li2020revisit}, etc. In \cite{Nesterov_book}, it is argued that the key behind constructing optimal methods lies in the accumulation of some global information on the objective function. The mathematical objects which enable for capturing the relevant topological information on the function that is to be minimized are the estimating sequences. Typically, they consist of a pair of sequences, that simultaneously allow for parsing global information around the iterates, as well as for measuring the convergence rate of the minimization process. Despite their remarkable properties, estimating sequences exhibit the issue that there is no unique or systematic approach for constructing them. As we will see in the sequel, making the adequate choice of the estimating functions that comprise the estimating sequences can significantly impact the practical performance of the resulting algorithm.  

An interesting framework for the study and analysis of various methods which are constructed based on the estimating sequences has been presented in \cite{baes2009estimate}. A unified analysis of accelerated first-order methods has appeared then in \cite{Tseng}. An existing estimating sequence method that can directly solve \eqref{opt_prob} is the Accelerated Multistep Gradient Scheme (AMGS) \cite{Nesterov_2007}. The method is proven to enjoy the accelerated rate of convergence $\mathcal{O}(\frac{1}{k^2})$. Despite its notable theoretical and practical performance as measured by the number of iterations carried through until convergence, the method suffers the drawback that it requires two projection-like operations per iteration. This results in an increase of the computational burden, which (in the case of large-scale problems) is also reflected in an increase of the runtime of the method. This problem has been solved by the development of the Fast Iterative Shrinkage-Thresholding Algorithm (FISTA) \cite{FISTA}. The method also enjoys the accelerated convergence rate of $\mathcal{O}(\frac{1}{k^2})$, while at the same time requiring only one projection-like operation per iteration. Similarly to \cite{Nesterov_83}, FISTA does not explicitly utilize the machinery of estimating sequences. However, as has been demonstrated in \cite{Iulian_2}, by properly selecting the estimating functions it is possible to establish links between FISTA and estimating sequence methods. 

Comparing AMGS to its counterpart devised for minimizing smooth convex functions \cite{Nesterov_book}, we can see that they were constructed using different types of estimating sequences. Nevertheless, they both enjoy the theoretical accelerated rate of convergence. However, although the theoretical analysis is unsuccessful in capturing the performance differences between these methods, their practical performance can vary quite significantly when they are tested on real-world problems and datasets. Moreover, from the simulations that we have performed for the cases of differentiable convex functions, we have noticed that the practical performance of FGM is better than AMGS and FISTA. Therefore, it is natural to ask whether the estimating sequences framework introduced in \cite{Nesterov_book} can be extended to the framework of composite objective functions, and evaluate the performance of the resulting method on problems of the form of \eqref{opt_prob}. In this work, we answer this question affirmatively, and show that, by constructing the appropriate estimating functions, it is possible to devise very efficient accelerated first-order methods. More specifically, the main contributions of the article are as follows.
\begin{itemize}
	\item From the theoretical perspective, we show how the gradient mapping technique \cite{Nemirovski_Yudin} can be coupled with the estimating sequences framework to produce a class of composite objective estimating-sequence techniques (COMET). Unlike AMGS, the resulting algorithms require only one projection-like operation per iteration.
	\item We introduce a tighter lower bound on the objective function than the one obtained from the Taylor series expansion of a convex function, and present a new structure for the estimating functions, which we call the composite estimating functions hereafter. We then show how they can be used to efficiently parse information around all the iterates, as well as measure the convergence rate of the minimization process. 
	\item Motivated by the recent improvement of the theoretical analysis of FGM for smooth convex functions \cite{dosti2020generalizing}, we show that by utilizing the proposed estimating sequences it is possible to achieve a better theoretical convergence rate even when the Lipschitz constant is not known and needs to be estimated. To compute the Lipschitz constant, we utilize a backtracking line search strategy. Moreover, we establish that the performance of the method is robust to the initialization of the parameters of the algorithm. 
	\item Through extensive simulations for various typical signal processing problems with composite structure, we show that the proposed method yields a better performance than the existing benchmarks. Furthermore, we also show the robustness of the selected instances of COMET with respect to the imperfect knowledge of the strong convexity parameter and the Lipschitz constant. To demonstrate the robustness, as well as the reliability of our proposed method, we test its performance on real-world datasets. 
\end{itemize}

The article is organized as follows. In Section~\ref{prel}, we introduce the key assumptions of the paper, as well as some of the main equations that are used in developing our method. Moreover, therein we also present the notation used throughout the paper. In Section~\ref{3}, we introduce the proposed estimating sequences for composite objectives and devise COMET based on them. In Section~\ref{ca}, we formally establish the convergence of COMET and derive the convergence rate for the minimization process. Then, in Section~\ref{ns}, we illustrate the practical performance of our proposed method in solving several optimization problems and show that it outperforms the existing benchmarks. Lastly, in Section~\ref{conc}, we present the conclusions. 

\section{Preliminaries and notation}
\label{prel}
We assume that the objective function is bounded below, i.e., \eqref{opt_prob} has a solution. The key assumption that we make is that the function and gradient computations have approximately the same complexity. In this framework, the necessary oracle functions are the function evaluators, $f(x)$, $g(x)$, gradient evaluator $\nabla f(x)$ and proximal evaluator $\text{prox}_{\tau g}(x)$. 

To simplify our analysis, let us relocate the strong convexity of $g(x)$ within the objective function in \eqref{opt_prob}. Let $x_0 \in \mathcal{R}^n$ and consider the following.
\begin{align}
	F(x) \! &= \! \left( \! f(x) + \frac{\tau \mu_g}{2} ||x - x_0||^2 \right) \! + \! \tau \! \left(g(x) - \frac{\mu_g}{2} ||x - x_0||^2\right) = \hat{f} (x) + \tau \hat{g} (x). \label{FFF}
\end{align}
The resulting function $\hat{f} (x)$ has a Lipschitz constant $L_{\hat{f}} = L_f + \tau \mu_g$ and strong convexity parameter $\mu_{\hat{f}} = \mu_f + \tau \mu_g$. On the other hand, the function $\hat{g} (x)$ has a strong convexity parameter $\mu_{\hat{g}} = 0$. 

Next, we recall that it is possible to construct upper and lower bounds for the smooth and strongly convex function $ \hat{f} (x)$ by using the following relations. 
\begin{align}
	\hat{f}(x) \leq \hat{f}(y) + \nabla \hat{f}(y)^T (x - y) + \frac{L_{\hat{f}}}{2} ||y - x||^2,  \label{upper_bound} \\
	\hat{f}(x) \geq \hat{f}(y) + \nabla \hat{f}(y)^T (x - y) + \frac{\mu_{\hat{f}}}{2} ||y - x||^2,  \label{lower_bound}
\end{align}
for all points $y \in \mathcal{R}^n$. Similarly, we can construct the following lower bound for the non-smooth term
\begin{align}
	\hat{g}(x) \geq \hat{g}(y) + s (y)^T (x - y),  \label{lower_bound_g}
\end{align}
where $s (y)$ is a subgradient of the function $\hat{g}(y)$. Moreover, for all $y \in \mathcal{R}^n$ and $L \geq L_{\hat{f}}$, we also define the following
\begin{align}
	\label{m_l}
	m_{L} (y;x) \triangleq \hat{f}(y) + \nabla \hat{f}(y)^T(x - y) + \frac{L}{2} ||x - y||^2 + \tau \hat{g}(x).
\end{align}
From the upper bound on the function established in \eqref{upper_bound}, we can write
\begin{align}
	\label{LLL_bound}
	m_{L}(y;x) \geq F(x), \forall x,y \in \mathcal{R}^n.
\end{align} 
At this point, we can introduce the composite gradient mapping as 
\begin{align}
	\label{T_L(y)}
	T_{L}(y) &\triangleq \arg \underset{x \in \mathcal{R}^n}{\min} \; m_{L} (y;x),
\end{align}
Lastly, we define the composite reduced gradient to be
\begin{align}
	\label{r_l}
	r_{L} (y) \triangleq L \left(y - T_{L}(y)\right).
\end{align}
Notice that when $\tau = 0$, we can write $\hat{f} (x) = f(x)$, which follows from \eqref{FFF}. Moreover, we have from \eqref{T_L(y)} that $T_L(y) = y - \frac{\nabla \hat{f}(y)}{L}$. Substituting this result into the definition given in \eqref{r_l}, yields $r_L (y) = \nabla F(y) = \nabla f(y)$, i.e., the composite reduced gradient becomes the gradient of the objective function. Moreover, from the first-order optimality conditions for \eqref{T_L(y)}, we can write 
\begin{align}
	\nonumber
	\nabla m_L(y;T_L(y))^T (x - T_L(y)) &\geq 0, \\ \label{r(y)}
	\left(\nabla f(y) + L(T_L (y) - y) + \tau s_L(y) \right)^T(x - T_L(y)) &\geq 0,
\end{align}
where $s_L (y) \in \partial F(T_L(y))$ is a subgradient belonging to the subdifferential of $F(T_L(y))$, whose value depends on the point $y$. Equating the first bracket of \eqref{r(y)} to $0$, as well as recalling definition \eqref{r_l}, we obtain the following relation useful for computing the value of the composite reduced gradient
\begin{align}
	\label{reduced_grad}
	r_L(y) = L(y - T_L(y)) = \nabla f(y) + \tau s_L(y).
\end{align}

We conclude the section by presenting a tighter lower bound on the objective function. 
\begin{theorem}
	\label{thm2}
	Let $F(x)$ be a composition of an $L_{\hat{f}}$-smooth and $\mu_{\hat{f}}$-strongly convex function $\hat{f}(x)$, and a simple convex function $\hat{g}(x)$. For $L \geq L_{\hat{f}}$, and $x,y \in \mathcal{R}^n$ we have 
	\begin{align}
		F(x) &\geq \hat{f}(T_L(y)) + \tau \hat{g} (T_L(y)) + r_L(y)^T \left(x - y\right) + \frac{\mu_{\hat{f}}}{2} ||x - y||^2 + \frac{1}{2L} ||r_L(y)||^2. \label{13}
	\end{align}
\end{theorem} 
\begin{proof}
	See Appendix \ref{A}. 
\end{proof}

\section{COMET}
\label{3}
In this section, we introduce our proposed method. We start by introducing the estimating sequences, and then show why these objects are useful. We also introduce a pair of estimating functions and show how they can be computed recursively. Then, based on the proposed construction of estimating functions, we derive COMET. 

We begin by defining the estimating sequences as
\begin{definition}
	\label{def__1}
	The sequences $\{\phi_{k}\}_{k = 0}^\infty$, and $\{\lambda_{k}\}_{k = 0}^\infty$, $\lambda_{k} \geq 0$, are called estimating sequences of the function $F(\cdot)$, if $\lambda_{k} \rightarrow 0$ as $k \rightarrow \infty$, and $\forall x \in \mathcal{R}^n$, $\forall k \geq 0$ we have
	\begin{equation}
		\label{def_1}
		\phi_{k} (x) \leq\lambda_{k} \phi_{0}(x) + (1 - \lambda_{k}) F(x).
	\end{equation}
\end{definition} 
Next, we proceed to proving that the estimating sequences allow for measuring the convergence rate to optimality.
\begin{lemma}
	\label{SFGM_lemma_1}
	If for some sequence of points $\{x_k\}_{k = 0}^\infty$ we have $F(x_k) \leq \phi_{k}^* \! \triangleq \! \underset{x \in {\mathcal{R}^n}}{ \min } \phi_{k} (x)$, then $F(x_k) - F(x^*) \leq \lambda_{k} \left[ \phi_{0}(x^*) - F(x^*) \right]$, where $x^* = \arg \underset{x \in \mathcal{R}^n}{\min} F(x)$.
\end{lemma}
\begin{proof}
	See Appendix \ref{B}.
\end{proof}

At this point, we are ready to show how the estimating functions can be defined recursively. 
\begin{lemma}
	\label{SFGM_lemma_2}
	Assume that there exists a sequence $\{\alpha_k\}_{k = 0}^\infty$, where $\alpha_{k} \in (0, 1)$ $\forall k$, $\sum_{k = 0}^{\infty} \alpha_{k} = \infty$, and an arbitrary sequence $\{y_k\}_{k = 0}^\infty$. Furthermore, let $\lambda_{0}$ = 1 and assume that the estimates $L_k$ of the Lipschitz constant $L_{\hat{f}}$ are selected in a way that inequality \eqref{upper_bound} is satisfied for all the iterates $x_k$ and $y_k$. Then, the sequences $\{\phi_{k}\}_{k = 0}^\infty$ and $\{\lambda_{k}\}_{k = 0}^\infty$, which are defined recursively as
	\begin{align}
		\label{lambda_recursive}
		\lambda_{k+1} &= (1 - \alpha_k) \lambda_{k}, \\ 
		\phi_{k+1} (x) \!  &= \! (1 \! - \! \alpha_k) \phi_{k} (x) \! + \! \alpha_{k} \! \! \left( \! \! F\left(T_{L_k} \! (y_k) \right) \! + \! \frac{1}{2L_k} ||r_{L_k} \! (y_k)||^2 \! \! \right) \nonumber \\ &+ \alpha_k \left( r_{L_k}(y_k)^T (x - y_k)  + \frac{\mu_f}{2} ||x - y_k||^2\right), \label{phi_k+1_SFGM}
	\end{align}
	are estimating sequences.
\end{lemma}
\begin{proof}
	See Appendix \ref{C}.  
\end{proof}
Let us now provide a comparison between the results obtained in Lemmas~\ref{SFGM_lemma_1}~and~\ref{SFGM_lemma_2} with their counterpart devised for the simpler case of minimizing smooth convex functions presented in \cite{Nesterov_book}. First, we can see from Lemma~\ref{SFGM_lemma_1} that the convergence rate of the minimization process depends entirely on the rate at which $\lambda_k \rightarrow 0$. Moreover, the result hints that for problem \eqref{opt_prob} we should expect a similar convergence rate as in the simpler case of minimizing a differentiable convex function. Then, in Lemma~\ref{SFGM_lemma_2}, we have shown how to form the estimating functions. From \eqref{phi_k+1_SFGM} we can see that we are utilizing a tighter lower bound than the one used for deriving FGM for the smooth strongly convex case.\footnote{To see why that is the case, recall that when $F(x)$ is smooth and convex function, the composite reduced gradient becomes just the gradient of the function.} Furthermore, we note that the cost function is evaluated at specific points in its domain, which are produced by the composite gradient mapping. Lastly, we can observe that we need to utilize the subgradient of the non-smooth objective function to construct the sequence of estimating functions $\{ \phi_k \}_{k=0}^\infty$. 

At this point, we note that no particular structure for the functions in the sequence $\{ \phi_k \}_{k=0}^\infty$ has been proposed yet. Inspired by the analysis for FGM in the setup of smooth convex functions \cite{Nesterov_book}, in the sequel we let
\begin{align}
	\label{phi}
	\phi_k(x) = \phi_0^* + \frac{\gamma_k}{2} ||x - v_k||^2,~ \forall k = 1, 2, \ldots,
\end{align}
which we call hereafter as scanning function. 
Nevertheless, we stress that this selection is not unique. As a matter of fact, different choices of the canonical structure for the function $\phi_k (x)$ can lead to entirely different algorithms, see for example \cite{pmlr-v75-zhang18a,ahn2020nesterov}. Now, the following Lemma is in order.
\begin{lemma}
	\label{SFGM_lemma_3}
	Let $\phi_{0} (x) = \phi_{0}^* + \frac{\gamma_{0}}{2} ||x - v_{0}||^2$. Then, the process defined in Lemma \ref{SFGM_lemma_2} preserves the canonical form of the scanning function presented in \eqref{phi}, where the sequences $\{\gamma_{k}\}_{k=0}^\infty$, $\{v_{k}\}_{k=0}^\infty$ and $\{\phi_{k}^*\}_{k=0}^\infty$ can be computed as follows
	\begin{align}
		\label{gamma_expr}
		\gamma_{k+1} &= (1-\alpha_k)\gamma_{k} + \alpha_k \mu_{\hat{f}}, \\ 
		\label{v_value}
		v_{k+1} \!  &= \! \frac{1}{\gamma_{k+1}}\left((1 \! - \! \alpha_k)\gamma_k v_{k} \! + \! \alpha_k \left(\! \mu_{\hat{f}} y_k \! - \! L_k \! \left(y_k \! - \! T_{L_k} \! \! \left( y_k \right) \! \right) \! \! \right) \! \!  \right) \! , \\
		\phi_{k+1}^*  &= (1-\alpha_k) \phi_k^* + \! \alpha_{k} \! \left( \! \! F\left(T_{L_k} (y_k) \right) \! + \! \frac{1}{2L_k} ||r_{L_k} (y_k)||^2 \! \! \right) - \! \frac{L_k^2 \alpha_k^2}{2 \gamma_{k+1}} ||y_k \! - \! T_{L_k}(y_k)||^2 \! \nonumber \\ &+ \! \frac{\mu_{\hat{f}} \alpha_k \gamma_k (1 \! - \! \alpha_k)}{2 \gamma_{k+1}} ||y_k \! - \! v_k||^2 + \frac{L_k \alpha_k \gamma_k (1-\alpha_k)}{\gamma_{k+1}} (y_k - v_k)^T(y_k - T_{L_k}(x_k)). \label{psi_{k+1}^*}
	\end{align}
\end{lemma} 
\begin{proof}
	See Appendix \ref{D}. 
\end{proof}
Comparing the result obtained in Lemma \ref{SFGM_lemma_3} with its counterpart constructed for minimizing smooth objective functions \cite[Lemma 2.2.3]{Nesterov_book}, we can see that the recursion for computing the elements in the sequences $\{v_{k}\}_{k=0}^\infty$ and $\{\phi_{k}^*\}_{k=0}^\infty$ has changed. It now reflects both the different lower bound on the objective function, as well as the reduced composite gradient, which were utilized for constructing the estimating functions that we now more precisely and informatively call as {\it composite estimating functions}. 

Now, we can proceed to constructing the algorithm via induction. First, let $\phi_0^* = F(x_0)$. Next, assume that for some iteration $k$, we have: $\phi_k^* \geq F(x_k)$. To conclude the induction argument, we need to establish that $\phi_{k+1}^* \geq F(x_{k+1})$. Utilizing our assumption for iteration $k$ in \eqref{psi_{k+1}^*}, we obtain
\begin{align}
	\phi_{k+1}^* \! &\geq \! (1 \! - \! \alpha_k) F(x_k) \! + \! \alpha_{k} \! \left( \! \! F\left(T_{L_k} (y_k) \right) \! + \! \frac{1}{2{L_k}} ||r_{L_k} \! (y_k)||^2 \! \! \right) - \! \frac{{L_k}^2 \alpha_k^2}{2 \gamma_{k+1}} ||y_k \! - \! T_{L_k} \! (y_k)||^2 \nonumber \\ &+ \! \frac{\mu_{\hat{f}} \alpha_k \gamma_k (1 \! - \! \alpha_k)}{2 \gamma_{k+1}} ||y_k \! - \! v_k||^2 + \frac{{L_k} \alpha_k \gamma_k (1-\alpha_k)}{\gamma_{k+1}} (v_k - y_k)^T(y_k - T_{L_k}(y_k)). \label{36}
\end{align}
Then, substituting the bound obtained in Theorem~\ref{thm2}, as well as \eqref{r_l} into \eqref{36}, we obtain 
	\begin{align}
		\phi_{k+1}^* \!  &\geq \! (1 \! - \! \alpha_k) \! \left( \! F(T_{L_k}(y_k)) \! + \! r_{L_k} \! (y_k)^T \! \left(x_k \! - \! y_k\right) \! + \! \frac{\mu}{2} ||x_k \! - \! y_k||^2 \! + \! \frac{1}{2L_k} ||r_{L_k} \! (y_k)||^2\right) \nonumber \\ & +\! \alpha_{k} \! \left( \! \! F\left(T_{L_k} (y_k) \right) \! + \! \frac{1}{2L_k} ||r_{L_k} (y_k)||^2 \! \! \right) - \! \frac{ \alpha_k^2}{2 \gamma_{k+1}} ||r_{L_k}(y_k)||^2 \! + \! \frac{\mu \alpha_k \gamma_k (1-\alpha_k}{2 \gamma_{k+1}} ||y_k - v_k||^2 \nonumber \\ &+ \frac{\alpha_k \gamma_k (1-\alpha_k)}{\gamma_{k+1}} r_{L_k}(y_k)^T (v_k - y_k). \label{37}
	\end{align}
	
Making some algebraic manipulations and factoring in \eqref{37}, we obtain
\begin{align}
	\phi_{k+1}^*  &\geq F(T_{L_k}(y_k)) + \left(\frac{1}{2{L_k}} -  \frac{\alpha_k^2}{2 \gamma_{k+1}}\right)||r_{L_k} (y_k)||^2 \nonumber \\ &+ (1-\alpha_k) r_{L_k}(y_k)^T \left(x_k - y_k + \frac{ \alpha_k \gamma_k }{\gamma_{k+1}} (v_k - y_k) \right). \label{38}
\end{align}

We are now ready to obtain a relation for the unknown terms in the sequences $\{\alpha_k\}_{k=0}^\infty$ and $\{y\}_{k=0}^\infty$. Observe that in \eqref{38} we can obtain the update rule for the terms in the sequence $\{\alpha_k\}_{k=0}^\infty$ as
\begin{align}
	\label{alpha_k_intuition}
	\alpha_k = \sqrt{\frac{\gamma_{k+1}}{L_k}}.
\end{align}
Utilizing the recursion for $\gamma _{k+1}$ given by \eqref{gamma_expr}, and solving the resulting quadratic equation yields
\begin{align}
	\label{alpha_k_SFGM}
	\alpha_{k} = \frac{\sqrt{\left(\mu_{\hat{f}} - \gamma_{k}\right)^2 + 4{L_k} \gamma_{k}}}{2{L_k}}.
\end{align}
Making the aforementioned selection for $\alpha_k$, we can now write \eqref{38} as
\begin{align}
		\phi_{k+1}^* &\geq F(T_{L_k}(y_k)) + (1-\alpha_k) r_{L_k}(y_k)^T \left(x_k - y_k + \frac{ \alpha_k \gamma_k }{\gamma_{k+1}} (v_k - y_k) \right). \nonumber
\end{align}
Therefore, we can obtain the update rule for the term $y_k$ by setting
\begin{align}
	x_k - y_k + \frac{ \alpha_k \gamma_k }{\gamma_{k+1}} (v_k - y_k) = 0.
\end{align}
This results in
\begin{align}
	\label{y_k}
	y_k &=  \frac{\gamma_{k+1} x_k + \alpha_k \gamma_{k} v_{k}}{\gamma_{k+1} + \alpha_k \gamma_{k}}. 
\end{align}
To establish that $\phi_{k+1} \geq F (x_{k+1})$, it suffices to let $x_{k+1} = T_{L_k} (y_k)$. 

Lastly, let us discuss another major difference between our proposed method and its counterpart derived for minimizing differentiable convex functions \cite{Nesterov_book}, which is the fact that our analysis allows for the line search adaptation.\footnote{Note that several backtracking strategies have already been proposed in the literature (see for instance \cite{FISTA, Tseng, Nesterov_2007})} The goal of our proposed line-search strategy is to select the smallest constant $L_k$ such that \eqref{upper_bound} is satisfied $\forall k = 0, 1, \ldots$. To progress faster towards $x^*$ in the initial iterations, we want to initialize $L_0 \in ]0, L_{\hat{f}}[$, and then slightly increase the value of the estimate of the Lipschitz constant across the iterations. 
However, since the true value of $L_{\hat{f}}$ is not known, this approach cannot be used. Therefore, it would be more preferable to select the line search strategy in a way that it ensures the robustness of the method with respect to the initialization of the estimate of the Lipschitz constant and ensure a dynamic update of the step size. Such a scheme would be of importance for many applications in signal processing (see \cite{Iulian_2} and the references therein). For this purpose, the following two parameters can be utilized: \textit{i)}~ a constant $\eta_u > 1$, which is used to increase the value of the estimate; \textit{ii)} a constant $\eta_d \in ]0,1[$, which is used to decrease the value of the estimate of the Lipschitz constant. Finally, the proposed method is summarized in Algorithm~\ref{FGM}.

\begin{algorithm}
	\caption{COMET}
	\label{FGM}
	\begin{algorithmic}[1]
		\STATE{\textbf{Input} $x_0 \in \mathcal{R}^n$, $L_0 > 0$, $\gamma_{0} \in [0, 3L_0 + \mu_{\hat{f}}]$, $\eta_u > 1$ \newline and $\eta_d \in ]0,1[$.}
		\STATE{\textbf{Set }$k = 0$, $i = 0$ and $v_0 = x_0$.}
		\WHILE{stopping criterion is not met}
		\STATE{$\hat{L}_i \leftarrow \eta_d L_k$}
		\WHILE{$F (\hat{x}_{i+1}) > m_{\hat{L}_i} (\hat{y}_i, \hat{x}_{i+1})$}
		\STATE{$
			\hat{\alpha}_{i} \leftarrow \frac{(\mu_{\hat{f}} - \gamma_{k}) + \sqrt{\! (\mu_{\hat{f}} - \gamma_{k})^2 + 4 \hat{L}_i \gamma_{k}}}{2 \hat{L}_i }$}
		\STATE{$ \label{133} \hat{\gamma}_{i+1} \leftarrow (1 - \hat{\alpha}_i) \gamma_{k} + \hat{\alpha}_i \mu_{\hat{f}} $}
		\STATE{$\hat{y}_i \leftarrow \frac{\hat{\gamma}_{i + 1} x_k + \hat{\alpha}_i \gamma_k v_k}{\hat{\gamma}_{i + 1} + \hat{\alpha}_i \gamma_k}$} 
		\STATE{ $\hat{x}_{i+1} \leftarrow \text{prox}_{\frac{1}{\hat{L}_i} \hat{g}} \left( \hat{y}_i - \frac{1}{\hat{L}_i} \nabla f(\hat{y}_i) \right)$ }
		\STATE{$\hat{v}_{i+1} \leftarrow \frac{1}{\hat{\gamma}_{i+1}}\left((1-\hat{\alpha}_i)\gamma_k v_{k} \! + \! \hat{\alpha}_i \left(\mu_{\hat{f}} \hat{y}_i - \hat{L}_i \left(\hat{y}_i \! - \! \hat{x}_{i+1}\right) \right)\right)$}
		\STATE{$\hat{L}_{i+1} \leftarrow \eta_u \hat{L}_i$}
		\STATE{$i \leftarrow i+1$}
		\ENDWHILE
		\STATE{$L_{k+1} \leftarrow \hat{L}_i$, $x_{k+1} \leftarrow \hat{x}_{i}$, $v_{k+1} \leftarrow \hat{v}_i$, $\alpha_k \leftarrow \hat{\alpha}_{i-1}$, $y_k \leftarrow \hat{y}_{i-1}$, $i \leftarrow 0$, $k \leftarrow k+1$}
		\ENDWHILE
		\STATE{\textbf{Output} $x_k$}
	\end{algorithmic}
\end{algorithm}

Let us now make a comparison between our proposed method and FGM (Constant Step Scheme I in \cite{Nesterov_book}). From lines 5 and 6 in Algorithm \ref{FGM} we can observe the similarities in updating the sequences $\{\alpha\}_{k=0}^\infty$ and $\{\gamma\}_{k=0}^\infty$. A difference can already be noticed in the update of the terms in the sequence $\{\gamma\}_{k=0}^\infty$, whose value becomes independent of $\mu_{\hat{f}}$. Then, a key difference between the methods is in the update of the iterates $x_k$. Due to the composite structure of the objective function of interest, the next iterate is computed by taking a proximal gradient step. Note that as long as the non-smooth term $g(x)$ has a simple structure, the proximal term can be computed efficiently. Another noticeable difference between the methods lies in the update of the terms in the sequence $\{v_k\}_{k=0}^\infty$, which now reflects the usage of the proposed subgradient. Lastly, we would like to point out that the parameter $\gamma_0$ can now be selected over a wider range of parameters than what is guaranteed by the existing convergence results for FGM established in \cite[Lemma~2.2.4]{Nesterov_book}. The rationale behind this result will become clear in the next section. 

Before we proceed to analyzing the convergence rate of the minimization process, let us evaluate the behavior of the estimate of the Lipschitz constant. Depending on the initialization of $L_0$, we obtain the following two scenarios: \textit{i)} If $L_0 \in ]0, L_{\hat{f}}[$, then from line 4 in Algorithm \ref{FGM}, we observe that the estimate of the Lipschitz constant at iteration $k$ increases only if $L_{k-1} \leq L_{\hat{f}}$. Therefore, we can write
\begin{align}
	\label{c1}
	L_0 \leq \hat{L}_i \leq L_k \leq \eta_u L_{\hat{f}}.
\end{align}
\textit{ii)} If $L_0 \geq L_{\hat{f}}$, then the condition in line~4 of Algorithm~\ref{FGM} is satisfied, and estimate of the Lipschitz constant cannot increase further. Therefore, we would have
\begin{align}
	\label{c2}
	L_k \leq \eta_d L_0.
\end{align}
Combining the bounds \eqref{c1} and \eqref{c2}, we can see that despite the initialization of $L_0$, we can always write
\begin{align}
	\label{L_bound}
	L_k \leq L_{\text{max}} \triangleq \text{max} \{\eta_d L_0, \eta_u L_{\hat{f}}\}.
\end{align}

\section{Convergence Analysis}
\label{ca}
Let us begin by noting that the result obtained in Lemma \ref{SFGM_lemma_1} suggests that the convergence rate of the minimization process will be the same as the rate at which $\lambda_k \rightarrow 0$. This is made more precise in the following theorem. 
\begin{theorem}
	\label{conv_analysis_t_1}
	If we let $\lambda_0 = 1$ and $\lambda_k = \prod_{i = 0}^{k-1} \left(1 - \alpha_i \right)$, Algorithm \ref{FGM} generates a sequence of points $\{x_k\}_{k = 0}^\infty$ such that
	\begin{align}
		F(x_k) - F^* &\leq \lambda_{k} \left[ F(x_0) - F(x^*) + \frac{\gamma_{0}}{2}||x_0 - x^*||^2 \right].
	\end{align} 
	\begin{proof} 
		See Appendix \ref{E}. 
	\end{proof}
\end{theorem}

We now recall that from Definition \ref{def__1}, we must have $\lambda_k \rightarrow 0$. Therefore, the result of Theorem \ref{conv_analysis_t_1} is sufficient to establish the fact that the sequence of iterates produced by our proposed algorithm converges to the optimal solution. The next step is to evaluate the rate of convergence of this process. Let us begin by characterizing the rate at which $\lambda_k \rightarrow 0$.
\begin{lemma}
	\label{conv_analysis_lemma_1}
	For all $k \geq 0$, Algorithm~\ref{FGM} guarantees that
	\begin{enumerate}
		\item If $\gamma_0 \in [0, \mu_{\hat{f}}[$, then 
		\begin{align}
			\label{FGM_conv_eq_66} 
			\lambda_{k} &\leq \frac{1.042 \mu_{\hat{f}} L_{\max}}{L_k^2 \left( e^{\frac{k + 1}{2} \sqrt{\frac{\mu_{\hat{f}}}{L_k}}} - e^{-\frac{k + 1}{2} \sqrt{\frac{\mu_{\hat{f}}}{L_k}}} \right)^2} \leq \frac{1.042 L_{\max}}{L_k (k+1)^2}.
		\end{align}
		\item If $\gamma_0 \in [\mu_{\hat{f}}, 3L_0 + \mu_{\hat{f}}]$, then
		\begin{align}
			\lambda_{k} \! &\leq \! \frac{4 \mu_{\hat{f}} L_{\max}}{L_k (\gamma_0 - \mu_{\hat{f}}) \! \left( \! e^{\frac{k + 1}{2} \sqrt{\frac{\mu_{\hat{f}}}{L_k}}} \! - e^{-\frac{k + 1}{2} \sqrt{\frac{\mu_{\hat{f}}}{L_k}}} \! \right)^2} \leq \frac{4L_{\max}}{\left[(\gamma_0 - \mu_{\hat{f}}) (k+1)\right]^2}.
			\label{FGM_conv_eq_66_second_interval} 
		\end{align}
	\end{enumerate}
\end{lemma}
\begin{proof}
	See Appendix~\ref{F}. 
\end{proof}
Comparing the results obtained in Lemma~\ref{SFGM_lemma_1} with the earlier results obtained in \cite[Lem\-ma~2.2.4]{Nesterov_book}, we can see two major differences. First, our proposed analysis establishes the convergence of the method even when the true value of the Lipschitz constant is not known.  Second, we can see that it is possible to establish the convergence of the method in minimizing objective functions with composite structure for a wider range of the initialization of the parameter $\gamma_0$. The importance of this result arises from the fact that the method exhibits the fastest theoretical and practical convergence when $\gamma_0 = 0$, which is not supported by the existing analysis of FGM. At the same time, the initialization $\gamma_0 = 0$ also provides more robustness to the imperfect knowledge of $\mu_{\hat{f}}$. 

From the result of Theorem~\ref{conv_analysis_t_1}, we can see that the convergence rate of the minimization process depends on the distance of the function value at $x_0$ to the function value at $x^*$. The following lemma yields an upper bound on the needed quantity.
\begin{lemma}
	\label{ll}
	Let $F(x)$ be a convex function with composite structure as shown in \eqref{opt_prob}. Moreover, let $T_L (y)$ and $r_L(y)$ be computed as given in \eqref{T_L(y)} and \eqref{reduced_grad}, respectively. Then, for any starting point $x_0$ in the domain of $F(x)$, we have
	\begin{align}
		F(x_0) - F(x^*) \leq \frac{L_0}{2} ||x_0 - x^*||^2. 
	\end{align}
	\begin{proof}
		See Appendix~\ref{G}. 
	\end{proof}
\end{lemma}

Combining the results of Lemmas \ref{conv_analysis_lemma_1} and \ref{ll} in Theorem \ref{conv_analysis_t_1}, we can immediately obtain the convergence rate for our method, as summarized in the following theorem.
\begin{theorem}
	\label{th3}
	Algorithm \ref{FGM} generates a sequence of points such that
	\begin{enumerate}
		\item If $\gamma_0 \in [0, \mu_{\hat{f}}[$, then
		\begin{align}
			F(x_k) - F(x_0) &\leq \frac{1.042 \mu_{\hat{f}} L_{\max} (L_0 + \gamma_0)||x_0 - x^*||^2}{2 L_k^2 \left( e^{\frac{k + 1}{2} \sqrt{\frac{\mu_{\hat{f}}}{L_k}}} - e^{-\frac{k + 1}{2} \sqrt{\frac{\mu_{\hat{f}}}{L_k}}} \right)^2} \leq \frac{1.042  L_{\max} (L_0 + \gamma_0) ||x_0 - x^*||^2}{2 L_k (k+1)^2}. \label{FGM_conv_eq_66} 
		\end{align}
		\item If $\gamma_0 \in [\mu_{\hat{f}}, 3 L_0 + \mu_{\hat{f}}]$, then
		\begin{align}
			F(x_k) \! - \! F(x_0) \! &\leq \! \frac{2 \mu_{\hat{f}}  L_{\max} (L_0 \! + \! \gamma_0)||x_0 \! - \! x^*||^2 }{L_k (\gamma_0 \! - \! \mu_{\hat{f}}) \! \left( \! e^{\! \frac{k + 1}{2} \! \sqrt{\frac{\mu}{L_k}}} \! - e^{\! -\frac{k \! + \! 1}{2} \! \sqrt{\frac{\mu}{L_k}}} \! \right)^2} \leq \frac{ 2L_{\max} (L_0 + \gamma_0)||x_0 - x^*||^2 }{\left[L_k (\gamma_0 - \mu_{\hat{f}})(k+1)\right]^2}. 	\label{FGM_conv_eq_66_second_interval} 
		\end{align} 
	\end{enumerate}
\end{theorem}
From the result of Theorem~\ref{th3} we can see that our proposed method is guaranteed to converge over a wider interval than its counterpart designed for minimizing smooth and strongly convex objectives. Notice that initializing $\gamma_0 = 0$ would guarantee the fastest convergence of the method. Such a result is important when considering many practical applications, wherein the true values of $\mu_{\hat{f}}$ and $L$ are often not known and should be estimated. The estimation error then propagates and affects the performance of the variants of COMET that would be initialized based on the estimates of the strong convexity and smoothness parameters. Another factor that impacts the rate of convergence of the minimization process is also the initialization of $L_0$. From \eqref{FGM_conv_eq_66} and \eqref{FGM_conv_eq_66_second_interval} we can see that the smaller $L_0$, the faster the convergence of the method. 

\section{Numerical study}
\label{ns}
In this section, we compare the numerical performance of the proposed method against the two state-of-the-art black-box methods, i.e., AMGS and FISTA, in solving several optimization problems, which arise often in many signal processing, statistics and data science applications. The selected loss functions are the quadratic and logistic loss functions, both with elastic net regularization. As we will see in the sequel, controlling the parameters of the elastic net regularizer allows for simulating extremely ill-conditioned examples. For the constructed examples, we show that COMET outperforms the selected benchmarks in terms of minimizing the number of iterations needed to achieve a certain tolerance level. To provide reliable results, we utilize both synthetic and real data, which are selected from the Library for Support Vector Machines (LIBSVM) \cite{LIBSVM}. To find the optimal solutions, we utilize CVX \cite{CVX}. 

In the first example, we illustrate the performance of three instances of COMET. More specifically, we consider the variant that in theory is expected to result in the fastest convergence, which is obtained when we initialize for $\gamma_{0} = 0$, and it is referred to as ``$\text{COMET}$, variant 1''. We also consider the variant that is expected to produce the slowest convergence, which happens when we initialize $\gamma_{0} = 3 L_0 + \mu_{\hat{f}}$, and it is labeled as ``$\text{COMET}$, variant 2''. Lastly, we also implement the instance of COMET that is obtained when $\gamma_0 = \mu_{\hat{f}}$, which is referred to as ``$\text{COMET}$, variant 3''. When comparing the performance of the methods under the condition where the Lipschitz constant is not known, for both AMGS and FISTA we utilize the line-search strategies presented in the respective works \cite{Nesterov_2007, FISTA}. We point out that throughout all the simulations the starting point $x_0$ is randomly selected and all algorithms are initialized in it.

\vspace{-3mm}
\subsection{Minimizing the quadratic loss function}
\label{linear_regression}
Consider one of the most popular problems in signal processing and statistics
\begin{equation}
	\begin{aligned}
		\label{mse}
		& \underset{x \in \mathcal{R}^n}{\text{minimize}}
		& &\frac{1}{2} \sum\limits_{i=1}^m (a_i^T x - y_i)^2 + \frac{\lambda}{2} ||x||^2 + \tau||x||_1,
	\end{aligned}
\end{equation}
where $|| \cdot||_1$ denotes the $l_1$ norm. The objective is to show that the theoretical gains of COMET, which are discussed in Section~\ref{ca}, are also reflected in the practical performance of the methods. Moreover, we analyze how the performance of the methods scales with the condition number of the problem. We also illustrate the practical benefits of utilizing the proposed line-search strategy. 

Let us first consider the simplest case, where the Lipschitz constant is assumed to be known. It allows for an objective assessment of the effectiveness of the methods in finding the optimal solution. For this example, we utilize synthetic data. We consider the diagonal matrix $A \in \mathcal{R}^{m \times m}$ and sample the elements $a_{ii}$ from the discrete set $\{10^0, 10^{-1}, 10^{-2}, \ldots, 10^{-\xi} \}$ uniformly at random. This choice of selecting $A$ ensures that $L = 1$ and $\mu_f = 10^{-\xi}$, which results in the condition number $10^{\xi}$. Then, we select the elements of the vector $y \in \mathcal{R}^m$ by uniformly drawing them from the box $[0, 1]^n$. Lastly, we note that in our computational experiments we set $m \in \{500, 1000, 1500, 2000\}$, $\xi \in \{3, 4, 7, 8\}$ and $\tau = \lambda \in \{ 10^{-3}, 10^{-4}, 10^{-7}, 10^{-8} \}$.

\begin{figure}
	\centering
	\begin{subfigure}[h]{0.49\columnwidth}
		\centering
		\includegraphics[width=1\columnwidth,height = \linewidth, trim={3cm 9.5cm 3cm 8cm}]{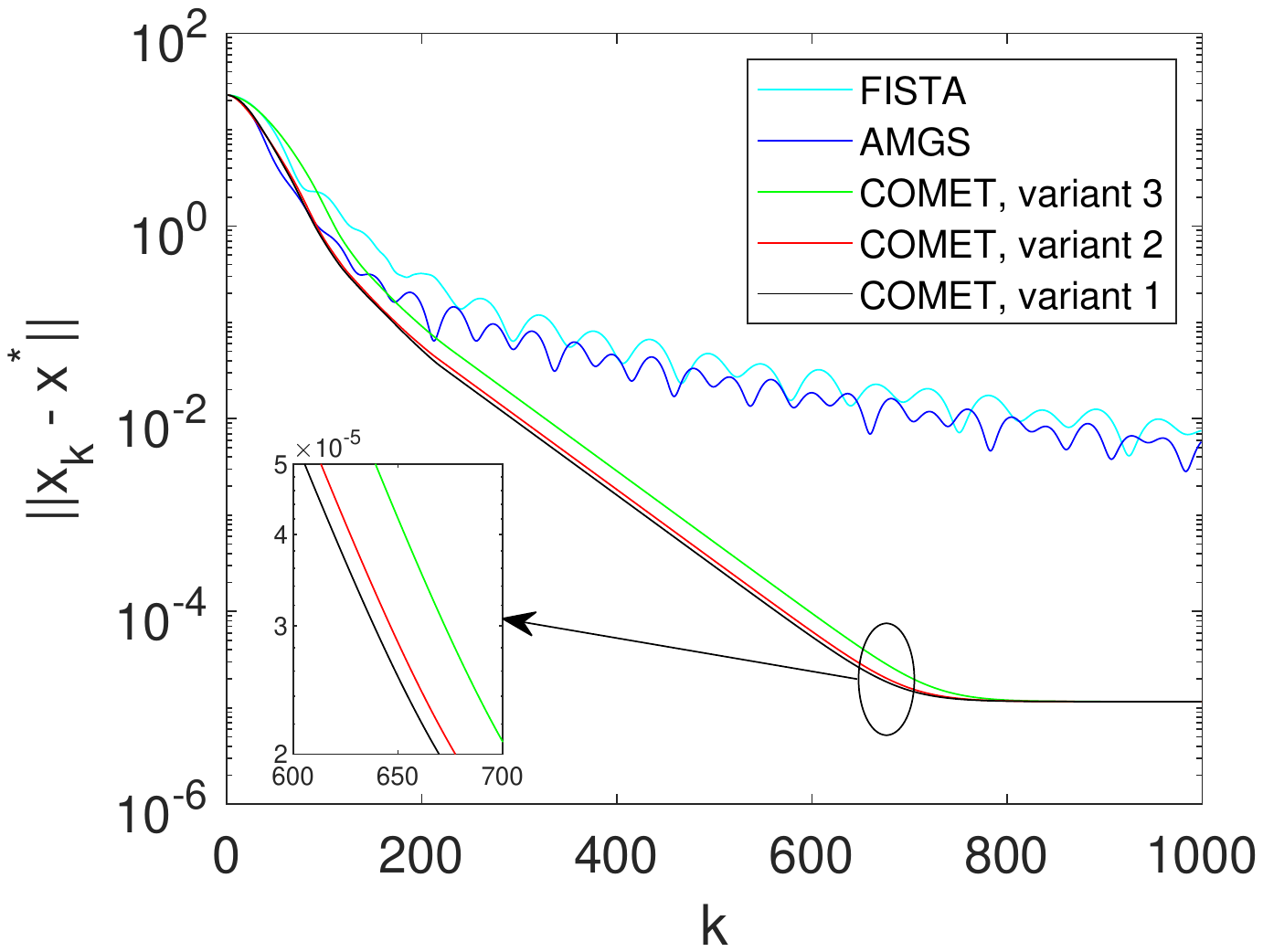} \label{num_sec_fig_1_sub_1} \vspace{-3mm}
		\caption{Decreasing the distance to $x^*$, $m = 500$, $\kappa = 10^3$ and $\lambda = \tau = 10^{-3}$.}
	\end{subfigure}
	\begin{subfigure}[h]{0.49\columnwidth}
		\centering
		\includegraphics[width=\columnwidth,height = \linewidth, trim={3cm 9.5cm 3cm 8cm}]{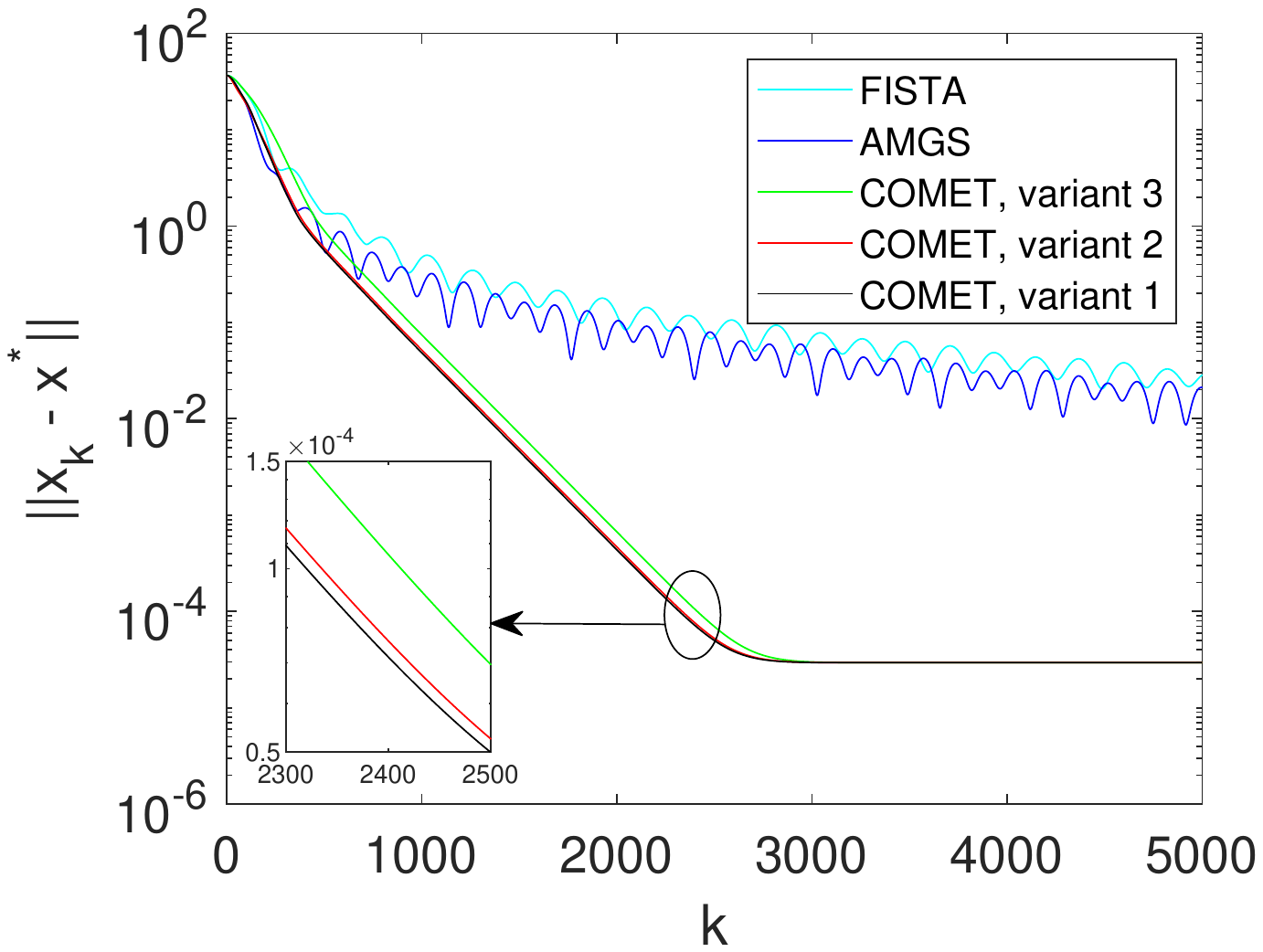} \label{num_sec_fig_1_sub_2} 
		\caption{Decreasing the distance to $x^*$, $m = 1000$, $\kappa = 10^4$ and $\lambda = \tau = 10^{-4}$.} 
	\end{subfigure}
	\begin{subfigure}[h]{0.49\columnwidth}
		\centering
		\includegraphics[width=1\columnwidth,height = \linewidth, trim={3cm 9.5cm 3cm 8cm}]{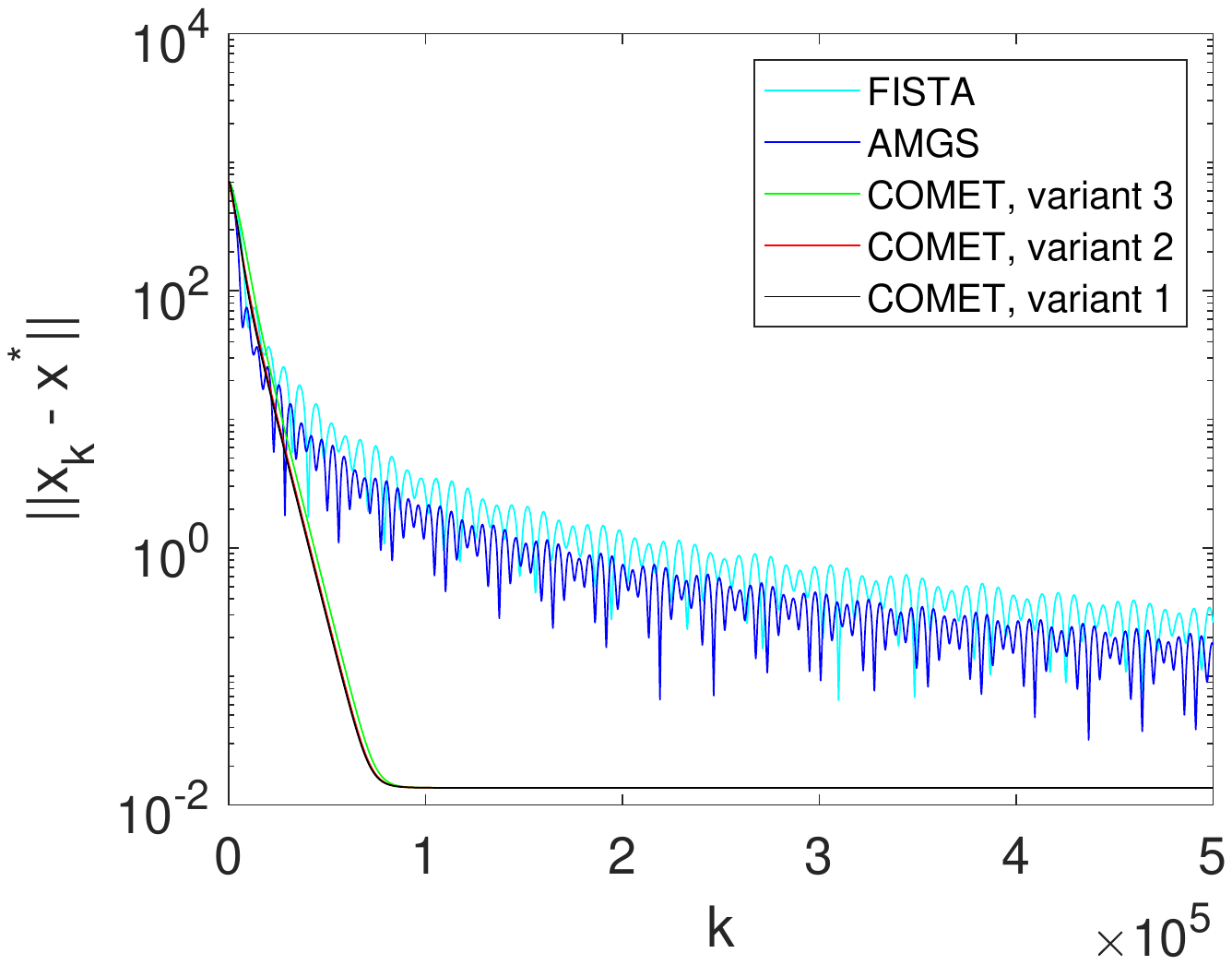} \label{num_sec_fig_1_sub_3}
		\caption{Decreasing the distance to $x^*$, $m = 1500$, $\kappa = 10^{7}$ and $\lambda = \tau = 10^{-7}$.}
	\end{subfigure}
	\begin{subfigure}[h]{0.49\columnwidth}
		\centering
		\includegraphics[width=\columnwidth,height = \linewidth, trim={3cm 9.5cm 3cm 8cm}]{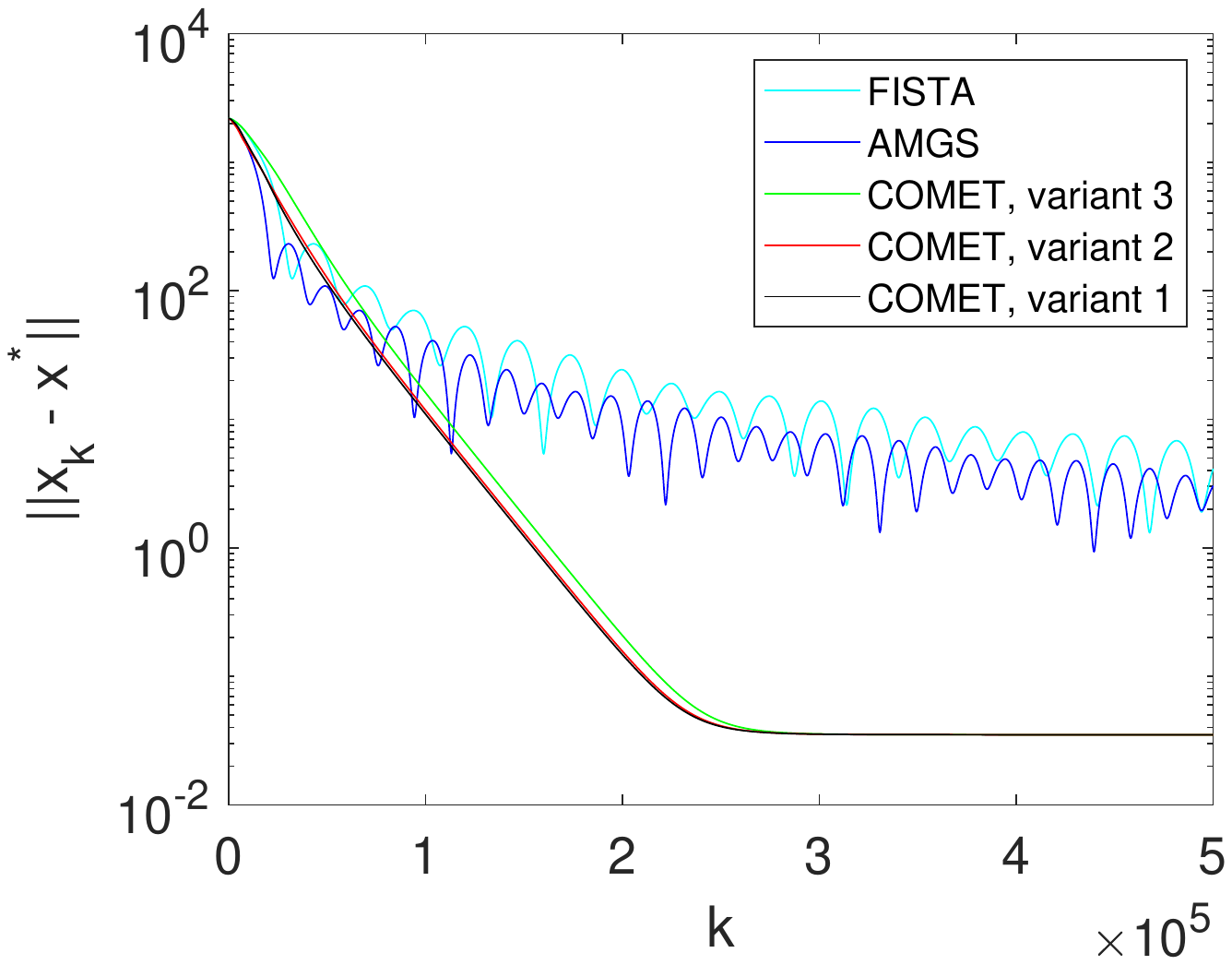} \label{num_sec_fig_1_sub_4}
		\caption{Decreasing the distance to $x^*$, $m = 2000$, $\kappa = 10^{8}$ and $\lambda = \tau = 10^{-8}$.} 
	\end{subfigure}
	\caption{Comparison between the efficiency of the algorithms tested in minimizing the quadratic loss function with elastic net regularizer on randomly generated data.}
	\label{num_sec_fig_1}
\end{figure}

From Fig.~\ref{num_sec_fig_1}, we can observe that the proposed method significantly outperforms all the existing benchmarks. First, notice that the larger the condition number of the problems becomes, the more iterations and consequently computations are required by the methods to obtain a good solution. Comparing between the methods, we can observe that all instances of COMET yield a better quality of the obtained solution, as measured by the distance to $x^*$. Moreover, we can clearly see that the iterates produced by COMET converge to $x^*$ in a much smaller number of iterations. Another important observation that can be made from the figure is that the proposed method exhibits better monotonic properties than both AMGS and FISTA. Comparing the performance of different variants of COMET, we can observe that their behavior is similar and the differences in performance are not too large. We can see that the variant that yields the best performance is the one obtained when $\gamma_0 = 0$, which is coherent with the theoretical results established in Section \ref{ca}. 

Next, we proceed to analyzing a more realistic scenario. We assume that the Lipschitz constant is not known, and needs to be estimated by using a line-search procedure. To demonstrate the robustness of the line-search strategy to be utilized in conjuction with COMET, we consider the following cases. \textit{i)} The Lipschitz constant is underestimated by a factor of $10$, i.e., $L_0 = 0.1L_{\hat{f}}$. \textit{ii)} The Lipschitz constant is overestimated by a factor of $10$, i.e., $L_0 = 10L_{\hat{f}}$. Moreover, we note that we selected $\eta_u = 2$ and $\eta_d = 0.9$, which were suggested in \cite{becker2011templates} because they ensure a good performance of the methods in many applications. Another parameter that is computationally expensive to be estimated in practice is the strong convexity parameter $\mu_{\hat{f}}$. To avoid an increase in computations, in all the following simulations we set the strong convexity parameter 
to be $0$. Lastly, we note that for all the examples that will be shown in the sequel, we utilize the datasets ``a1a'' and ``colon-cancer''. The former dataset has data matrix $A \in \mathcal{R}^{1605 \times 123}$, whereas the latter has $A \in \mathcal{R}^{62 \times 2000}$. 

\begin{figure}
	\centering
	\begin{subfigure}[h]{0.49\columnwidth}
		\centering
		\includegraphics[width=1\columnwidth,height = \linewidth, trim={3cm 9.5cm 3cm 8cm}]{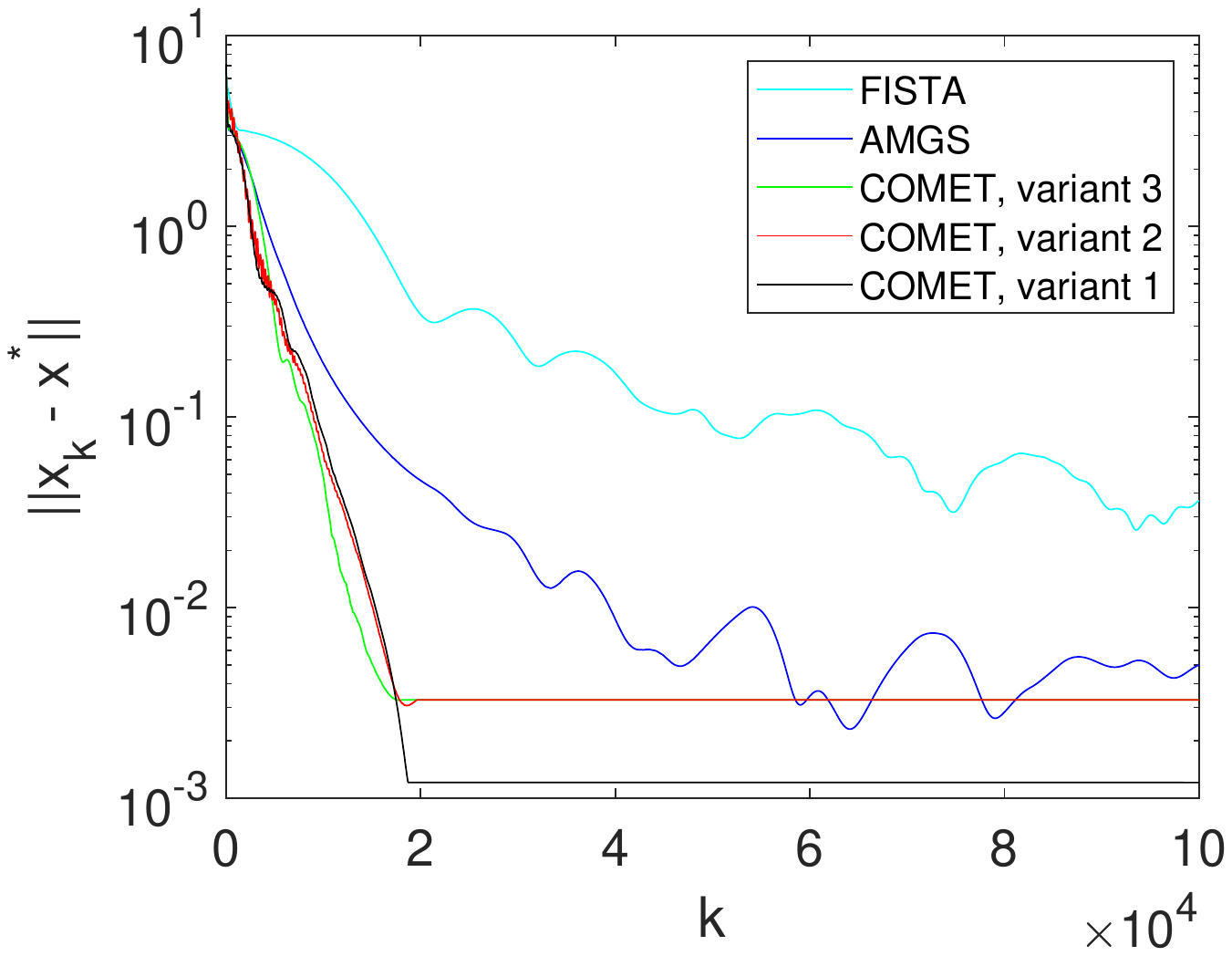} \label{num_sec_fig_2_sub_1} 
		\caption{Decreasing the distance to $x^*$ on ``a1a'' dataset, $L_0 = 0.1 L_{\text{``a1a''}}$ and $\lambda = \tau = 10^{-5}$.}
	\end{subfigure}
	\begin{subfigure}[h]{0.49\columnwidth}
		\centering
		\includegraphics[width=\columnwidth,height = \linewidth, trim={3cm 9.5cm 3cm 8cm}]{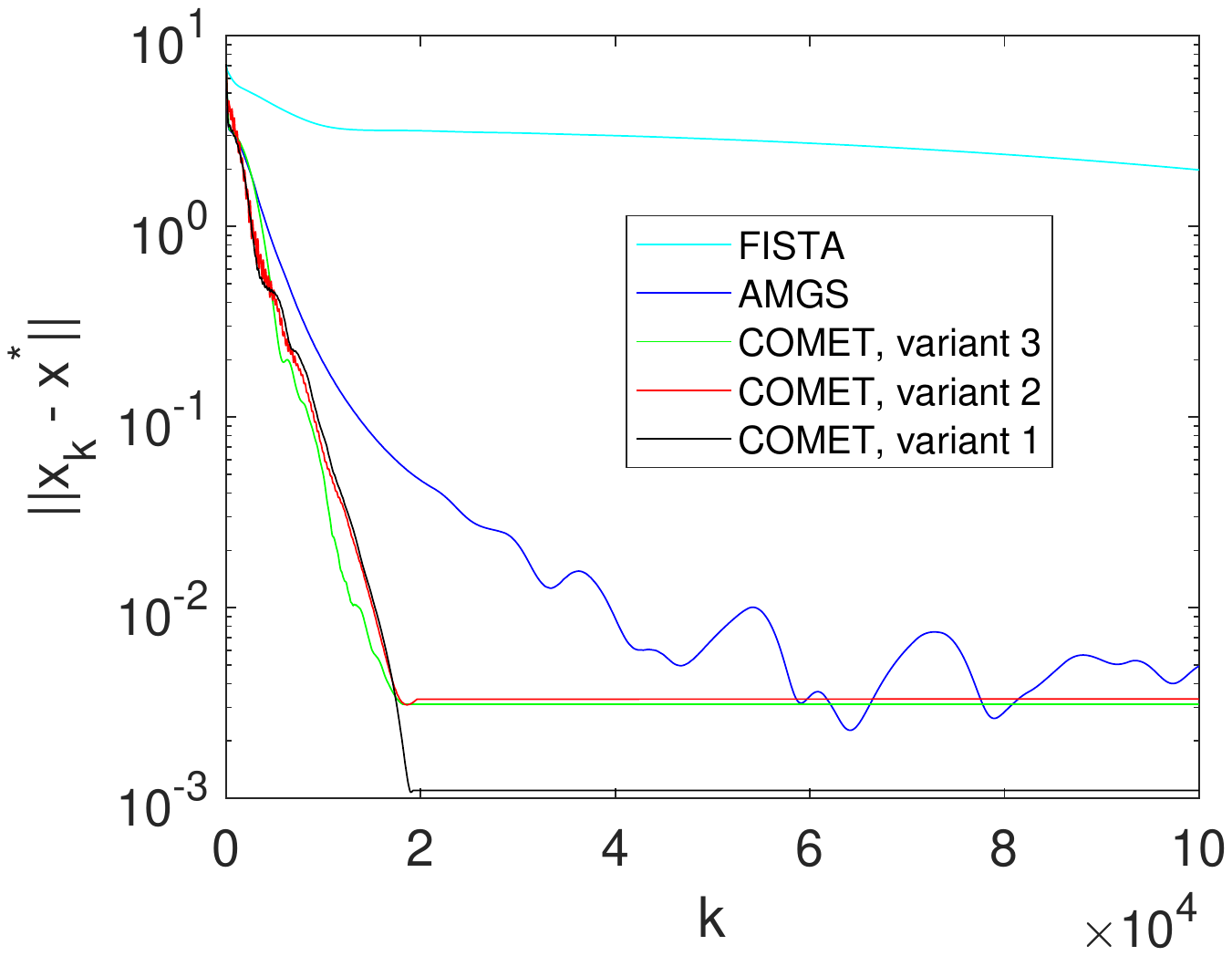} \label{num_sec_fig_2_sub_2} 
		\caption{Decreasing the distance to $x^*$ on ``a1a'' dataset, $L_0 =10 L_{\text{``a1a''}}$ and $\lambda = \tau = 10^{-5}$.}
	\end{subfigure}
	\begin{subfigure}[h]{0.49\columnwidth}
		\centering
		\includegraphics[width = \columnwidth,height = \linewidth, trim={3cm 9.5cm 3cm 8cm}]{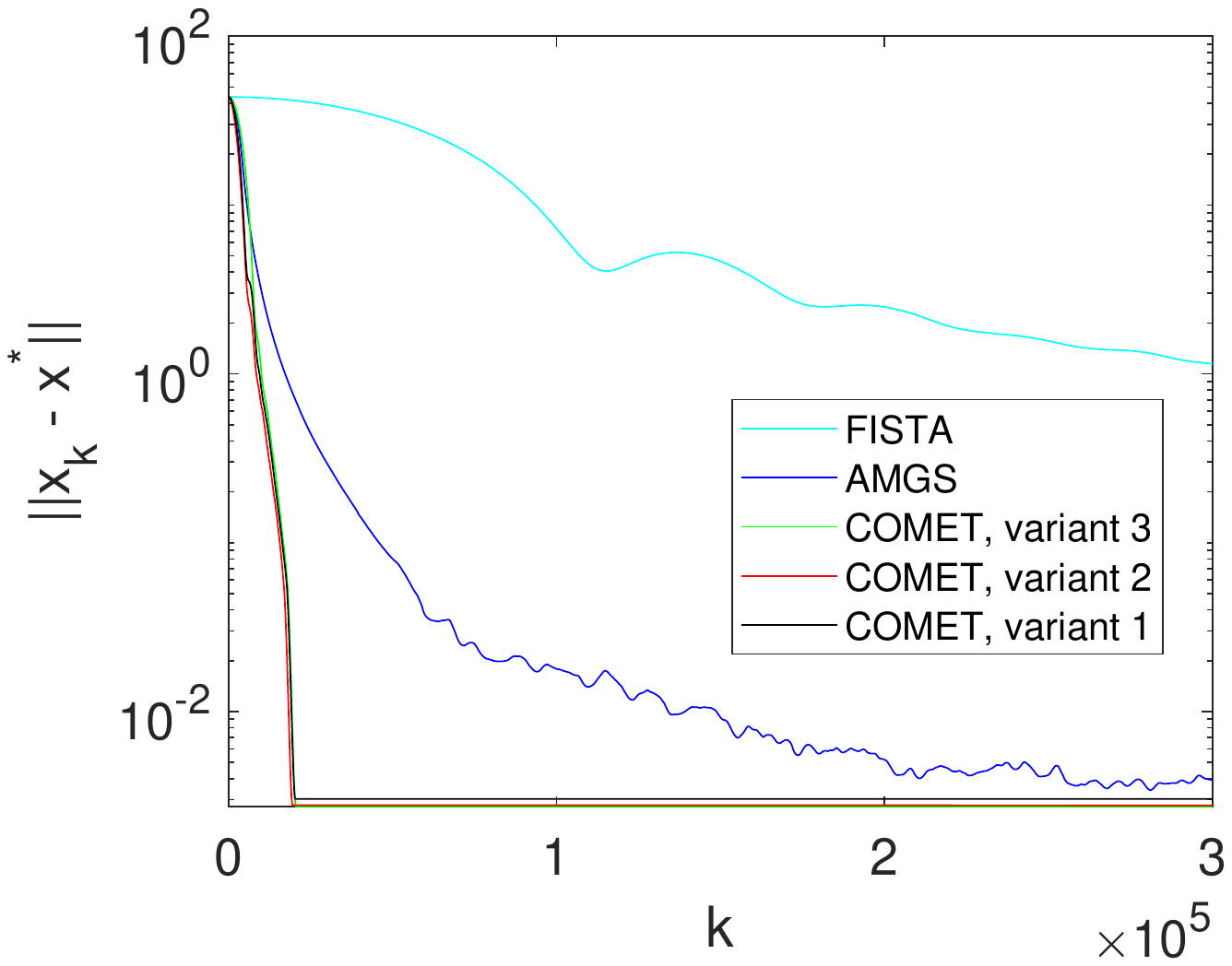} \label{num_sec_fig_2_sub_5} 
		\caption{Decreasing the distance to $x^*$ on ``colon-cancer'' dataset, $L_0 = 0.1 L_{\text{``colon-cancer''}}$ and $\lambda = \tau = 10^{-6}$.}
	\end{subfigure}
	\begin{subfigure}[h]{0.49\columnwidth}
		\centering
		\includegraphics[width=\columnwidth,height = \linewidth, trim={3cm 9.5cm 3cm 8cm}]{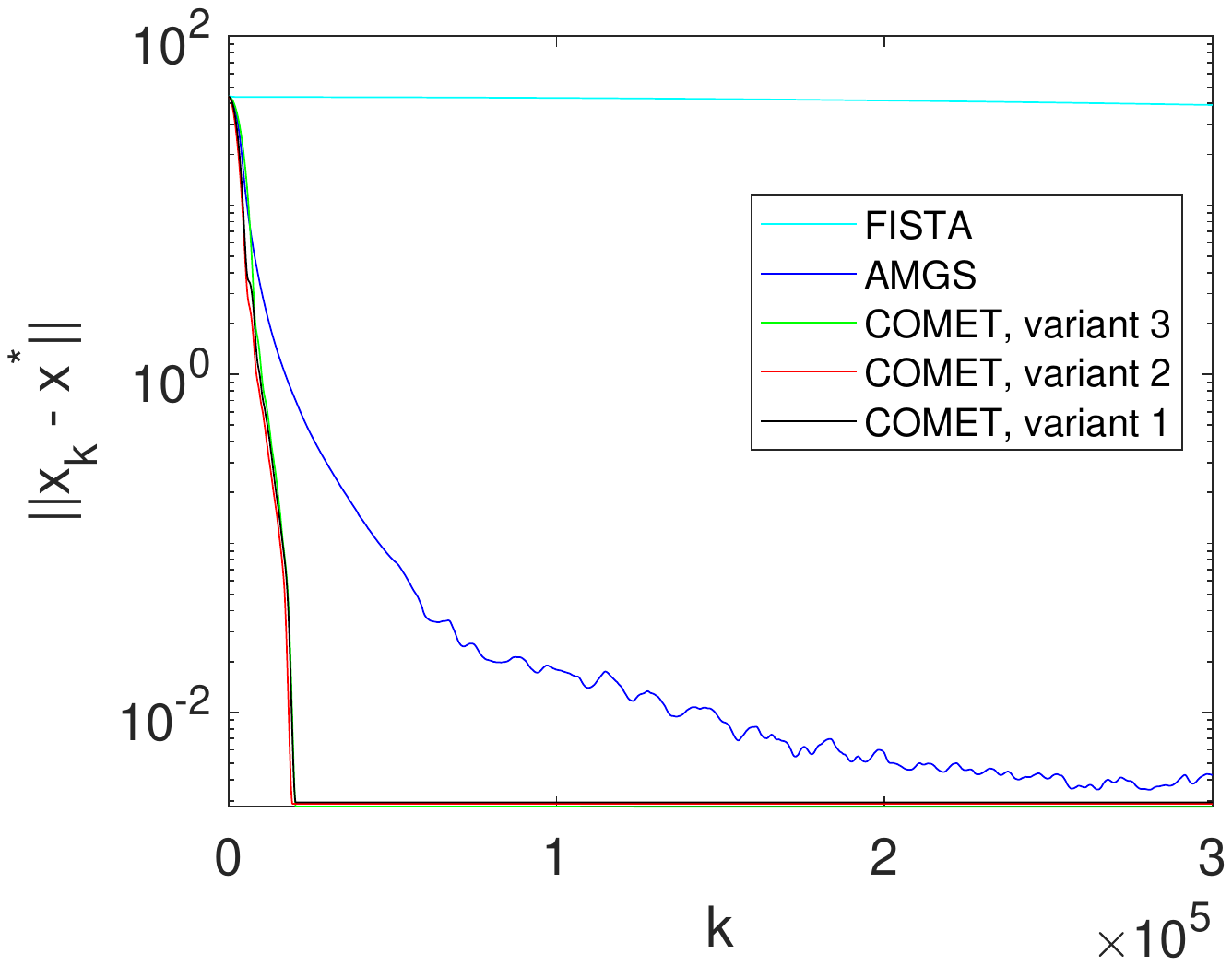} \label{num_sec_fig_2_sub_4} 
		\caption{Decreasing the distance to $x^*$ on ``colon-cancer'' dataset, $L_0 = 10 L_{\text{``colon-cancer''}}$ and $\lambda = \tau = 10^{-6}$.}
	\end{subfigure}
	\caption{Comparison between the efficiency and robustness with respect to the initialization of the Lipschitz constant of the algorithms tested in minimizing the quadratic loss function with elastic net regularizer on real data.} 
	\label{num_sec_fig_2}
\end{figure}

For the datasets that we are utilizing, the respective Lipschitz constants are $L_{\text{``a1a''}} = 10061$ and $L_{\text{``colon-cancer''}} = 1927.4$. Moreover, we let the regularizer term $\lambda = \tau \in \{10^{-5}, 10^{-6}\}$. Evidently, this selection of the regularizer terms guarantees a very large condition number for the problems that are being solved. The numerical results are presented in Fig.~\ref{num_sec_fig_2}, from which we can observe that all the instances of COMET significantly outperform the existing benchmarks. First, the final iterate produced by the first variant of $\text{COMET}$ is the closest to $x^*$. This is most visible from the numerical experiments conducted on the "a1a" dataset, which are depicted in Figs.~\ref{num_sec_fig_2}(a)~and~\ref{num_sec_fig_2}(b). Second, the iterates produced by the proposed COMET converge to $x^*$ by requiring a significantly smaller number of iterations, when compared to AMGS and FISTA. Third, the performance of FISTA largely depends on the initialization of the Lipschitz constant. On the other hand, we can observe that for both datasets, the performance of both AMGS and COMET remains unaffected by the value of $L_0$. We stress that COMET retains the robustness to $L_0$ at the lower computational cost of only one projection-like operation per iteration, whereas AMGS requires double of that. Lastly, comparing the performance between the selected variants of COMET, we can see that in practice their performance differences are minor. Nevertheless, our results shown in Figs.~\ref{num_sec_fig_2}(a)~and~\ref{num_sec_fig_2}(b) suggest that the version of COMET which is obtained when $\gamma_0 = 0$ yields a better performance. This becomes important particularly when considering practical applications, wherein the true values of $ \mu_{\hat{f}}$ and $L_{\hat{f}}$ are typically not known and their true values can only be estimated within some error bounds. From this perspective, we can conclude that the instance of COMET obtained by setting $\gamma_0 = 0$ enjoys both the faster convergence of the iterates and the robustness with respect to the imperfect knowledge of $\mu_{\hat{f}}$ and $L$.

\subsection{Minimizing the logistic loss function}
\label{Decreasing the norm of the gradient}
\begin{figure}
	\centering
	\begin{subfigure}[h]{0.49\columnwidth}
		\centering
		\includegraphics[width=1\columnwidth,height = \linewidth, trim={3cm 9.5cm 3cm 8cm}]{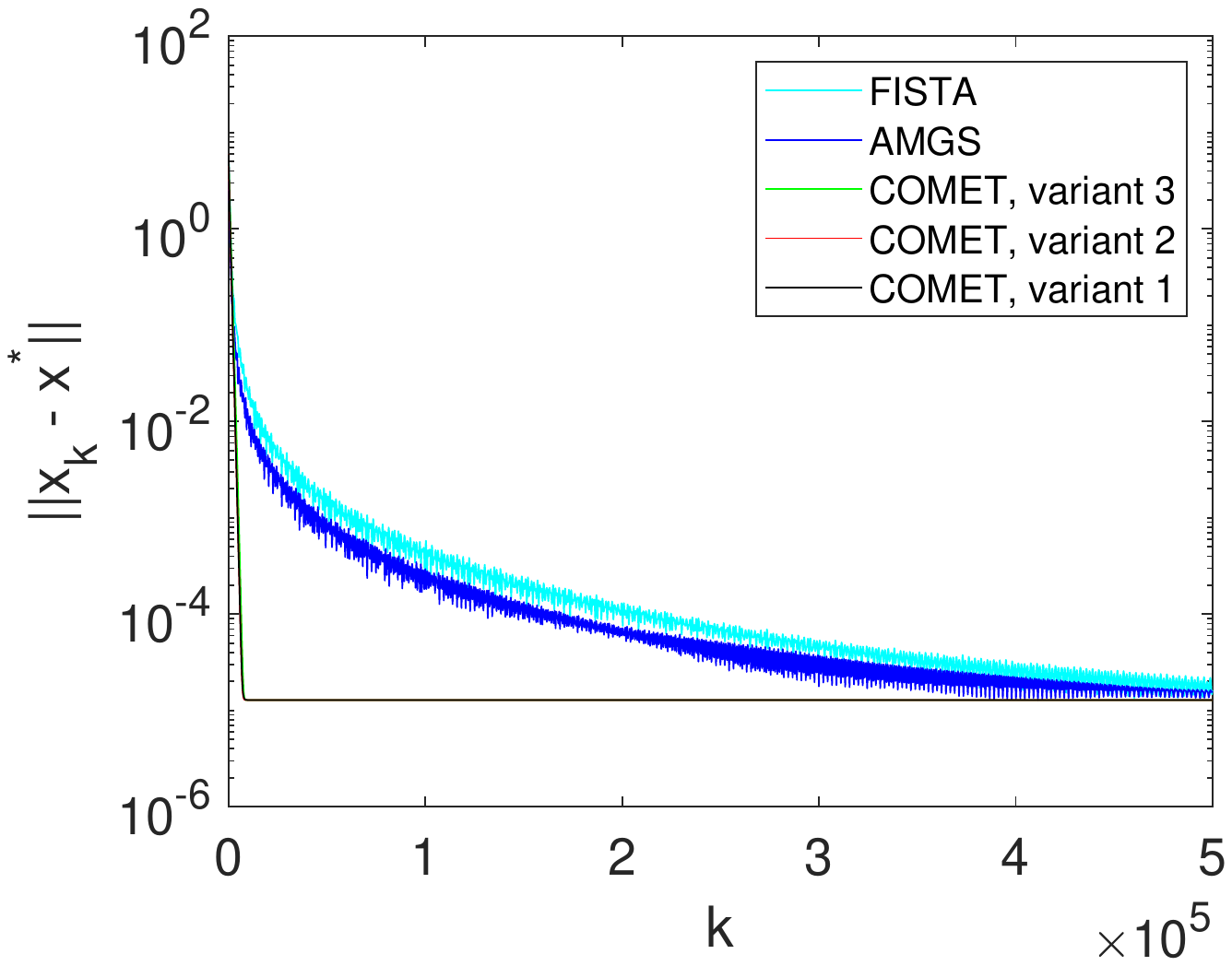} \label{num_sec_fig_3_sub_1} 
		\caption{Decreasing the distance to $x^*$ on ``triazine'' dataset, $\lambda = 10^{-4}$ and $\tau = 10^{-5}$.}
	\end{subfigure}
	\begin{subfigure}[h]{0.49\columnwidth}
		\centering
		\includegraphics[width=\columnwidth,height = \linewidth, trim={3cm 9.5cm 3cm 7.5cm}]{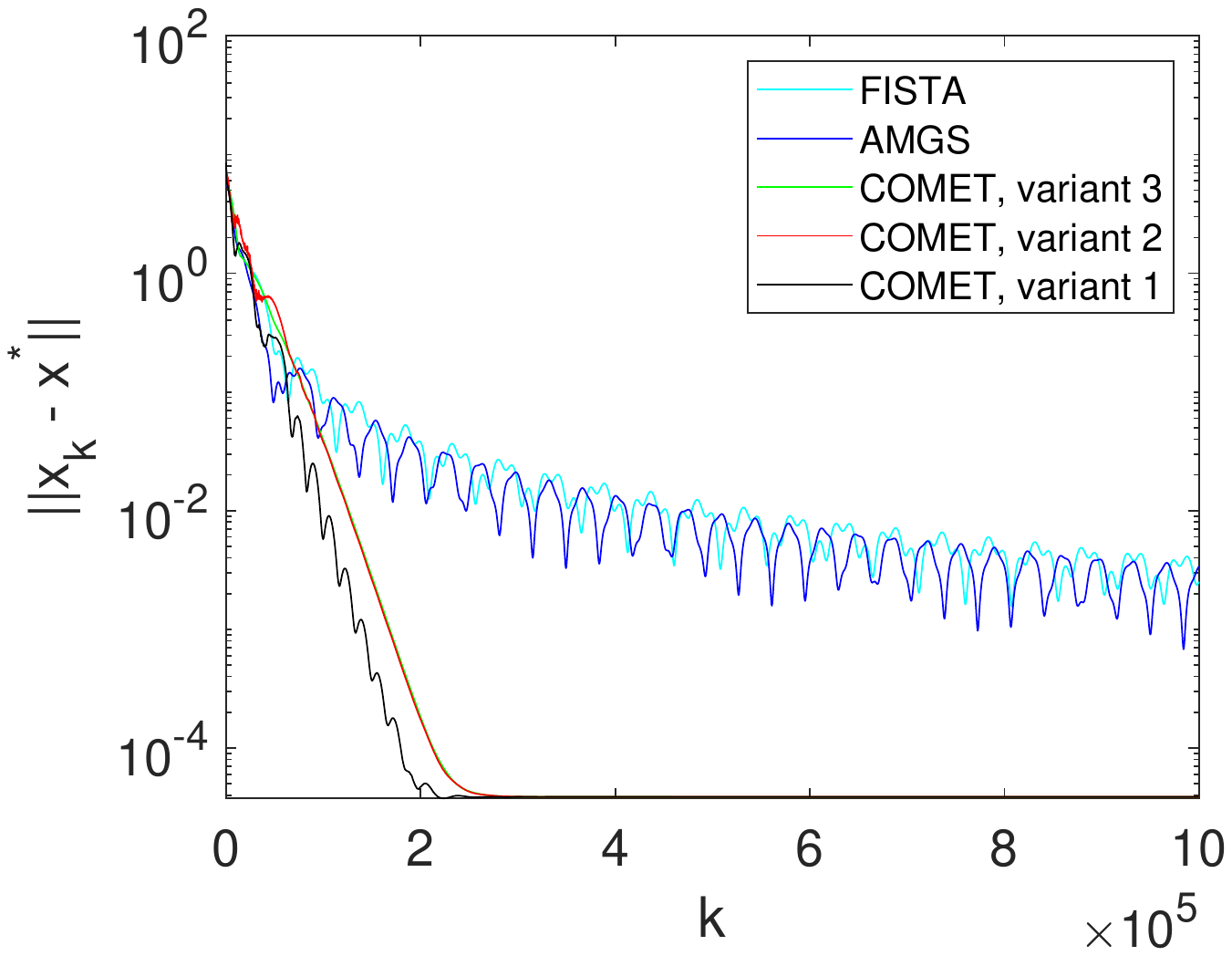} \label{num_sec_fig_3_sub_2} 
		\caption{Decreasing the distance to $x^*$ on ``triazine'' dataset, $\lambda = 10^{-7}$ and $\tau = 10^{-6}$.}
	\end{subfigure}
	\begin{subfigure}[h]{0.49\columnwidth}
		\centering
		\includegraphics[width = \columnwidth,height = \linewidth, trim={3cm 9.5cm 3cm 7.5cm}]{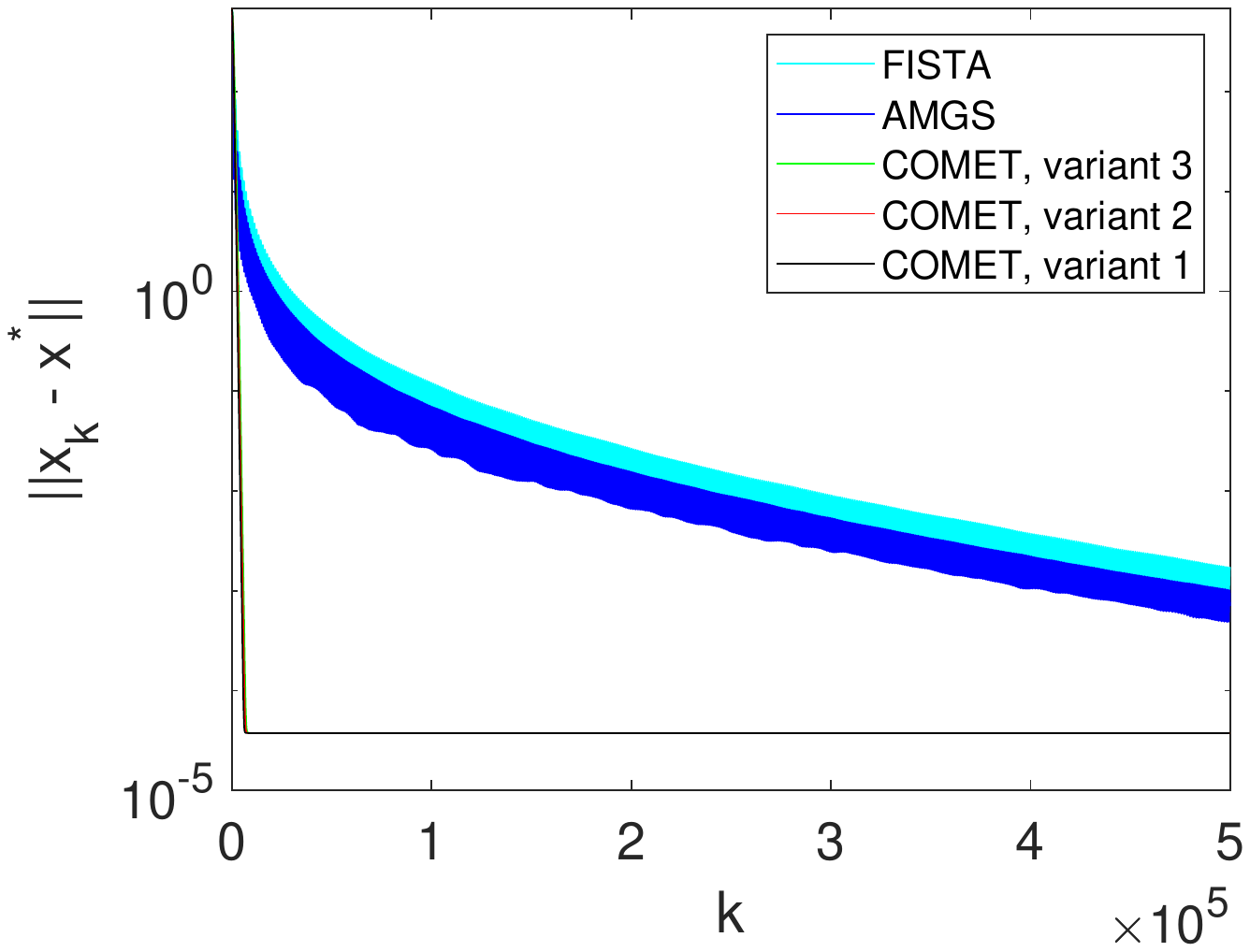} \label{num_sec_fig_3_sub_3} 
		\caption{Decreasing the distance to $x^*$ on a subset of ``rcv1.binary'' dataset, $\lambda = 10^{-6}$ and $\tau = 10^{-4}$.}
	\end{subfigure}
	\begin{subfigure}[h]{0.49\columnwidth}
		\centering
		\includegraphics[width=\columnwidth,height = \linewidth, trim={3cm 9.5cm 3cm 8cm}]{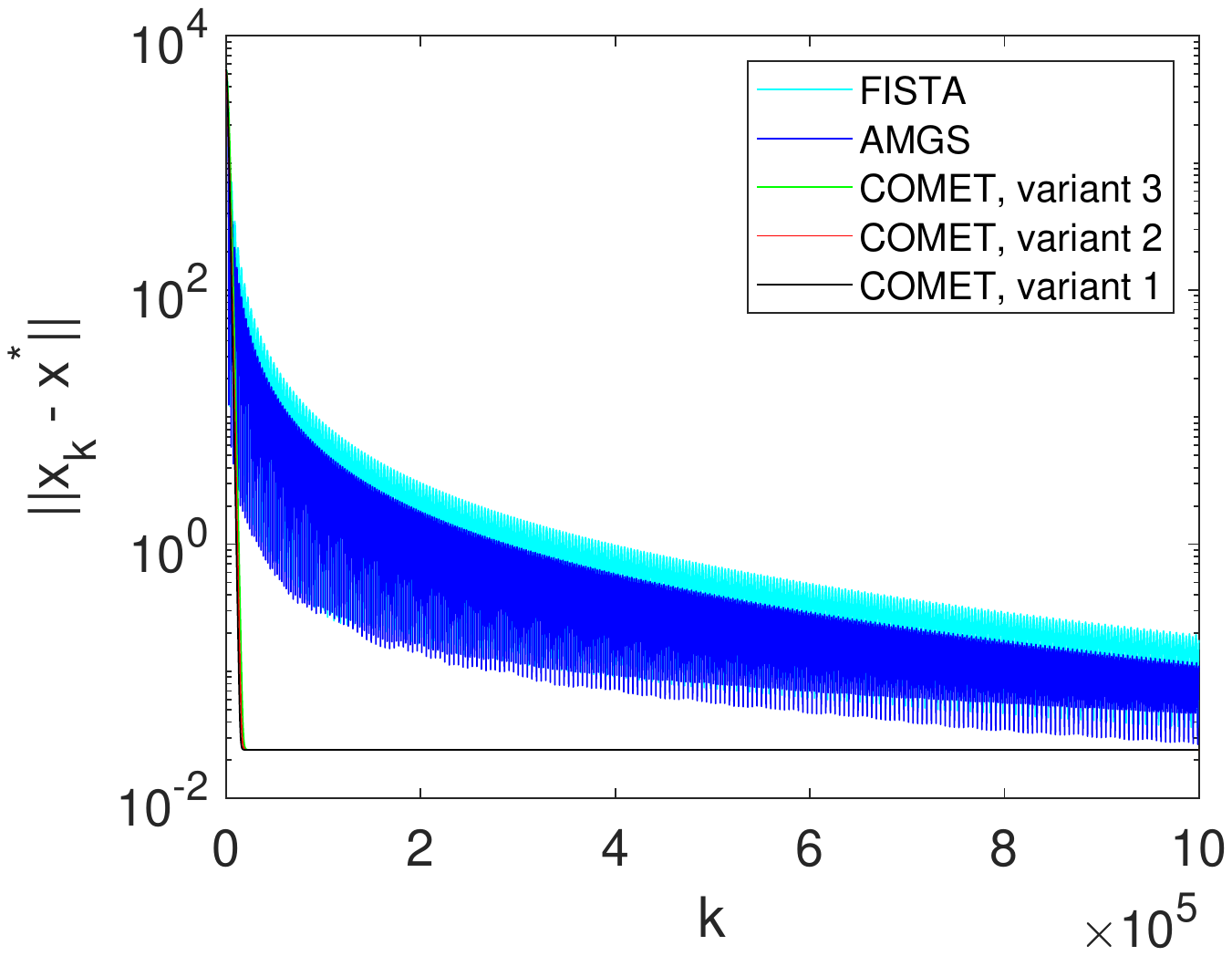} \label{num_sec_fig_3_sub_4} 
		\caption{Decreasing the distance to $x^*$ on a subset of ``rcv1.binary'' dataset, $\lambda = 10^{-5}$ and $\tau = 10^{-7}$.}	\end{subfigure}
	\caption{Comparison between the efficiency of the algorithms tested in minimizing the logistic loss function with elastic net regularizer on real data.}
	\label{num_sec_fig_3}
\end{figure}

To demonstrate the versatility of the proposed black-box method, let us now compare its performance to the selected benchmarks in minimizing a regularized logistic loss function with elastic net regularizer
\begin{equation}
	\begin{aligned}
		\label {num_eq_2}
		& \underset{x \in \mathcal{R}^n}{\text{minimize}}
		& &\frac{1}{m} \sum\limits_{i=1}^m \text{log} \left( 1 + \text{e}^{-b_i x a_i} \right) + \frac{\lambda}{2} ||x||^2 + \tau ||x||_1.  
	\end{aligned}
\end{equation} 
For this problem type, we diversify the utilized datasets and select ``triazine'', as well as a subset of ``rcv1.binary''. For the chosen datasets, we have $A_{\text{``triazine''}} \in \mathcal{R}^{186 \times 61}$ and $A_{\text{``rcv1.binary''}} \in \mathcal{R}^{1000 \times 2000}$. Moreover, from the results of Fig.~\ref{num_sec_fig_2}, we have observed that the performance of FISTA has been dependent on the initial estimate of the Lipschitz constant and has been overall worsened when $L_{\hat{f}}$ is unknown. Therefore, to provide the fairest comparison with respect to FISTA, for this set of examples we estimate the value of $L$ directly from the data. More specifically, we have $L_{\text{``triazine''}} = 25.15$ and $L_{\text{``rcv1.binary''}} = 1.13$. On the other hand, the strong convexity parameter of the data is set $\mu_{\hat{f}} = 0$. Lastly, we note that for this set of numerical experiments we consider the cases when $\lambda \neq \tau$. The results are reported in Fig. \ref{num_sec_fig_3}, wherein the specific values for $\lambda$ and $\tau$ are also presented. 

From Fig. \ref{num_sec_fig_3}, we can observe that for both datasets, COMET outperforms and exhibits better monotonic properties than AMGS or FISTA. Moreover, all variants of COMET require a much lower number of iterations to produce iterations which are closest to $x^*$. Lastly, for the selected problem type, the variant of COMET which is constructed when $\gamma_0 = 0$ yields the best practical performance, although the true value of $\mu_{\hat{f}}$ is not known.

\section{Conclusions and Discussion}
\label{conc}
We have considered the problem of constructing accelerated black-box first-order methods for solving optimization problems with composite structure by utilizing the estimating sequences framework. We have introduced a new class of estimating functions and shown that by exploiting their coupling with the gradient mapping technique it is possible to construct very efficient gradient-based methods, which we named COMET. Unlike the existing results on the convergence of FGM-type methods devised in \cite{Nesterov_book}, the novel convergence analysis established in this work allows for the adaptation of the step-size. Another major contribution which stemmed from the proposed convergence analysis is the fact that COMET is guaranteed to converge when $\gamma_0 \in [0, 3L + \mu_{\hat{f}}]$. The practical implication of these two observations is the fact that it is possible to construct efficient accelerated methods which are also robust to the imperfect knowledge of the smoothness and strong convexity parameters. Our theoretical findings are corroborated by extensive numerical experiments conducted on several problem types, wherein both synthetic and real-world data were utilized. 

{The results that were established in this work can be further developed in different directions. Particularly, it is interesting to investigate the possibilities of embedding the heavy-ball momentum into COMET. Another attractive research direction is the investigation of the possibility of coupling between the proposed framework and the inexact oracle framework, as well as the framework for constructing distributed proximal gradient methods. Lastly, we note that it is also interesting to investigate the possible extensions to designing accelerated algorithms for solving non-convex optimization problems.}

\appendices
\section{Proof of Theorem 1}
\label{A}
We start by showing that $m_L (y;x)$ is an $L$-strongly convex function in $x$. Notice that it is defined to be the sum of convex functions. Therefore, it is itself a convex function. Now, consider
\begin{align}
	\label{sfgm_eq_3}
	m_L(y;y) - m_L(y; T_L(y)) \geq \frac{L}{2} ||y - T_L(y)||^2. 
\end{align}
By the definition given in \eqref{T_L(y)}, $T_L(y)$ is the minimizer of $m_L(y;x)$ over all $x \in \mathcal{R}^n$. Therefore, we can conclude that $m_L (y; x)$ is a strongly convex function with strong convexity parameter $L$ \cite[Theorem 2.1.8]{Nesterov_book}.

Now, we can proceed to deriving the lower bound. From \eqref{lower_bound} and \eqref{lower_bound_g}, we can write
\begin{align}
	F(x) &\geq {\hat{f}}(y) + \tau {\hat{g}}(y) + \left( \nabla f(y) + \tau s_L (y) \right)^T \left(x - y\right) + \frac{\mu_{\hat{f}}}{2} ||x - y||^2 . \label{dfff}
\end{align}
Then, from the definition of $m_L (y,y)$ introduced in \eqref{m_l}, as well as \eqref{reduced_grad}, we can rewrite the right-hand side (RHS) of \eqref{dfff} as
\begin{align}
	&{\hat{f}}(y) + \tau {\hat{g}}(y) + \left(\nabla {\hat{f}}(y) + \tau s_L (y)\right)^T \! \left(x - y\right) + \frac{\mu_{\hat{f}}}{2} ||x - y||^2 \nonumber \\ &= m_L(y;y) + r_L(y)^T \left(x - y\right) + \frac{\mu_{\hat{f}}}{2} ||x - y||^2. \label{ttt}
\end{align}
Moreover, from \eqref{sfgm_eq_3}, we can write 
\begin{align}
	m_L(y;y) &+ r_L(y)^T \left(x - y\right) \geq m_L(y; T_L(y)) + \frac{L}{2} ||y - T_L(y)||^2 + r_L(y)^T \left(x - y\right). \nonumber
\end{align}
Utilizing the definition of the reduced composite gradient introduced in \eqref{r_l}, yields
\begin{align}
	&m_L(y; T_L \! (y)) \! + \! \frac{L}{2} ||y \! - \! T_L(y)||^2 \! + \! r_L(y)^T \! \!  \left(x \! - \! y\right) \! + \! \frac{\mu_{\hat{f}}}{2} ||x \! - \! y||^2 \nonumber \\ &= \! m_L(y; T_L \! (y)) \! + \! \frac{1}{2 L} ||r_L(y)||^2 \! + \! r_L \! (y)^T \! \! \left(x \! - \! y\right) \! + \! \frac{\mu_{\hat{f}}}{2} ||x \! - \! y||^2 \! \!. \label{tttt}
\end{align}
Finally, we note that taking a proximal gradient descent step on $f(x)$, which by assumption has Lipschitz continuous gradient, we can derive \eqref{13}. This completes the proof. 

\section{Proof of Lemma 1}
\label{B}
By the assumption of Lemma \ref{SFGM_lemma_1}, we have
\begin{align}
	\nonumber
	F(x_k) \leq \phi_{k}^* &= \underset{x \in {\mathcal{R}^n}}{\min} \phi_{k} (x) \stackrel{\eqref{def_1}}{\leq}  \!
	\underset{x \in {\mathcal{R}^n}}{\min} \left[\lambda_{k} \phi_{0}(x) + (1 - \lambda_{k}) F(x) \right] \leq  \lambda_{k} \phi_{0}(x^*) + (1 - \lambda_{k}) F(x^*).
\end{align}
Rearranging the terms yields the desired result. 

\section{Proof of Lemma 2}
\label{C}
We prove this lemma by induction. Let us begin by analyzing iteration $k = 0$. By assumption, we have $\lambda_0 = 1$. Utilizing \eqref{def__1}, we obtain $\phi_0(x) \leq \lambda_0 \phi_0(x) + \left( 1 - \lambda_0 \right) F(x) \equiv \Phi_0 (x)$. Then, we assume that \eqref{def__1} holds true at some iteration $k$. Therefore, we can write
\begin{align}
	\label{useful}
	\phi_k(x) - \left(1 - \lambda_k \right) F(x) \leq \lambda_k \phi_0 (x).
\end{align}

Substituting the bound obtained in Theorem \ref{thm2}, i.e., \eqref{13} in \eqref{phi_k+1_SFGM}, we obtain  
\begin{align}
	\label{crappy}
	\phi_{k+1} (x) &\leq (1 - \alpha_k) \phi_{k} (x) + \alpha_k F(x).
\end{align}
Then, adding and subtracting the same term to the RHS of \eqref{crappy}, we reach 
\begin{align}
	\phi_{k+1} (x) \! &\leq \! (1 \! - \! \alpha_k) \phi_{k} (x) \! + \! \alpha_k F(x) \! + \! (1 - \alpha_{k}) (1 - \lambda_{k}) F(x) - (1 - \alpha_{k}) (1 - \lambda_{k}) F(x) \nonumber \\ 
	&= (1 - \alpha_k) \left[\phi_{k} (x) - (1 - \lambda_{k}) F(x)\right] + \left(\alpha_k + (1 - \lambda_{k}) (1 - \alpha_k) \right) F(x). \label{3.8}
\end{align}
Using the bound obtained in \eqref{useful} in \eqref{3.8}, we have
\begin{align}
	\label{mmmm}
	\phi_{k+1} (x) &\leq (1 \! - \! \alpha_k) \lambda_{k} \phi_{0} (x) + (1 - \lambda_{k} + \alpha_{k} \lambda_{k}) F(x).
\end{align}     
Lastly, after utilizing \eqref{lambda_recursive}, the proof is concluded.

\section{Proof of Lemma 3}
\label{D}
Let us begin with establishing the first part of the proof through a mathematical induction argument. At iteration $k = 0$, we have $\nabla^2 \phi_0 (x) = \gamma_0 I$. Next, we assume that at some iteration $k$ it is true that $\nabla^2 \phi_k(x) = \gamma_k I$. Then, at iteration $k+1$ we can write
\begin{align}
	\label{t}
	\nabla^2 \phi_{k+1} (x) &\stackrel{\eqref{phi_k+1_SFGM}}{=} (1-\alpha_k) \gamma_{k}  I + \alpha_k \mu_{\hat{f}} I  \equiv \gamma_{k+1} I. 
\end{align}

Now, we proceed to establishing the proposed recurrent relations for updating the terms in the sequences $\{v_{k}\}_{k=0}^\infty$ and $\{\phi_{k}^*\}_{k=0}^\infty$. Substituting \eqref{phi} into \eqref{phi_k+1_SFGM}, and analyzing its first-order optimality conditions we obtain 
\begin{align}
		\gamma_{k+1} (x - v_{k+1})  &=  \gamma_k  (1  -  \alpha_k) (x  - v_{k})  + \alpha_k  \left(\mu_{\hat{f}} (x  - y_k) + \; r_{L_k}(y_k) \right). \label{SFGM_opt_cond}
\end{align} 
We can then reduce the terms that depend on $x$ by utilizing \eqref{gamma_expr} in \eqref{SFGM_opt_cond}. This results in
\begin{align}
	\label{v}
	\begin{split}
		- \gamma_{k+1} v_{k+1} &= - (1 - \alpha_k) \gamma_k v_{k} + \alpha_k \left(- \mu_{\hat{f}}  y_k + r_{L_k}(y_k)\right).
	\end{split}
\end{align} 
Then, substituting \eqref{r_l} in \eqref{v}, we obtain \eqref{v_value}. 

We will conclude this proof by establishing \eqref{psi_{k+1}^*}. Let us begin by substituting \eqref{phi} in \eqref{phi_k+1_SFGM}, now evaluated at the points $x = y_k$. This way we obtain 
	\begin{equation}
		\label{gen_conditions_thetas}
		\phi_{k+1}^* + \frac{\gamma_{k+1}}{2}||y_k - v_{k+1}||^2 = (1-\alpha_k)\left(\phi_k^* \!+\! \frac{\gamma_k}{2}||y_k - v_k||^2 \right) \!+\! \alpha_{k} \left( F\left(T_{L_k} (y_k) \right) + \frac{1}{2L_k} ||r_{L_k} (y_k)||^2 \right). 
	\end{equation}
Next, we proceed by utilizing \eqref{v_value} to compute the second term in the left hand side (LHS) of \eqref{gen_conditions_thetas}. Consider the following
\begin{align}
	v_{k+1}  - y_k &= \frac{1}{\gamma_{k+1}}((1-\alpha_k)\gamma_k v_{k} \! + \! \alpha_k \mu_{\hat{f}} y_k - \alpha_k L_k \left(y_k \! - \! T_{L_k} \left( y_k \right)\right) - \gamma_{k+1} y_k). \label{vv}
\end{align}
Then, utilizing \eqref{gamma_expr} in \eqref{vv}, we obtain
\begin{align}
	\label{vv}
	v_{k+1} \! - \! y_k \! &= \! \frac{1}{\gamma_{k+1}}((1 \! - \! \alpha_k)\gamma_k (v_{k} \! - \! y_k) \! - \! \alpha_k L_k \! \left(y_k \! - \! T_{L_k} \! \left( y_k \right)\right) \! . 
\end{align}
Taking $||\cdot||^2$ of both sides in \eqref{vv}, yields
\begin{equation}
	\label{more_eq}
	||y_k \!  - \!  v_{k+1}||^2 \!  =  \! \frac{||(1 \! - \! \alpha_k)\gamma_k (v_{k} \! - \! y_k) \! - \! \alpha_k L_k \! \left(y_k \! - \! T_{L_k} \!  \! \left( y_k \right)\right) \!  ||^2}{\gamma_{k+1}^2} \!  . 
\end{equation}
Finally, multiplying both sides of \eqref{more_eq} by $\frac{\gamma_{k+1}}{2}$ and expanding the RHS, we obtain
	\begin{align}
		\label{fgm_eq_1}
		\frac{\gamma_{k+1}}{2}||y_k - v_{k+1}||^2 \! &= \! \frac{(1 \! - \! \alpha_k)^2\gamma_k^2}{2 \gamma_{k+1}} ||v_{k} \! - \! y_k||^2 + \frac{\alpha_{k}^2 {L_k}^2}{2 \gamma_{k+1}} ||y_k - T_{L_k} \left( y_k \right) ||^2 \nonumber \\ &- \frac{2 {L_k} \alpha_k (1-\alpha_k)\gamma_{k}}{2 \gamma_{k+1}} (v_{k} - y_k)^T\nabla \left(y_k - T_{L_k} \left( y_k \right)\right) \!.
	\end{align}
Substituting \eqref{fgm_eq_1} in \eqref{gen_conditions_thetas}, and making some straightforward algebraic manipulations, we obtain \eqref{psi_{k+1}^*}.

\section{Proof of Theorem 2}
\label{E}
Set $\phi_0^* = f(x_0)$. Then, considering \eqref{phi} evaluated at iteration $k = 0$, we obtain $\phi_0(x_0) = f(x_0) + \frac{\gamma_0}{2} ||x_0 - v_0||^2$. In Algorithm \ref{FGM}, we initialize $v_0 = x_0$, which is sufficient to guarantee that $f(x_0) \leq \phi_0^*$ at step $k= 0$. Moreover, recall that we designed the update rules of the proposed method to guarantee that $f(x_k) \leq \phi_k^* $, $\forall k = 1,2, \ldots$. Therefore, the necessary conditions for the results proved in Lemma \ref{SFGM_lemma_1} to be applied are satisfied. 

\section{Proof of Lemma 4}
\label{F}
Let $\gamma_0 \in [0, 3L_0 + \mu_{\hat{f}}]$ and consider applying the result obtained in \eqref{gamma_expr} to the following
\begin{align}
	\label{k}
	\gamma_{k+1} - \mu_{\hat{f}} &= (1 - \alpha_k) \gamma_{k} + \alpha_k \mu_{\hat{f}} - \mu_{\hat{f}}. 
\end{align}
Then, utilizing the assumption that $\lambda_0 = 1$ in \eqref{k}, we can write 
\begin{align}
	\label{kkkk}
	\gamma_{k+1} - \mu_{\hat{f}} = (1 - \alpha_k) \lambda_0 \left[  \gamma_{k} - \mu_{\hat{f}} \right].
\end{align}
Using the recursivity of \eqref{gamma_expr} in relation \eqref{kkkk}, yields 
\begin{align}
	\label{FGM_conv_eq_1}
	\gamma_{k+1} - \mu_{\hat{f}} = \lambda_{k+1}\left[\gamma_{0} - \mu_{\hat{f}} \right].
\end{align}
Let us now exploit the connection between relations \eqref{lambda_recursive} and \eqref{alpha_k_intuition}, which can be linked through the term $\alpha_k$ as follows
\begin{align}
	\label{t}
	\alpha_k &= 1 - \frac{\lambda_{k+1}}{\lambda_k} = \sqrt{\frac{\gamma_{k+1}}{L_k}} = \sqrt{\frac{\mu_{\hat{f}}}{L_k} + \frac{\gamma_{k+1}- \mu_{\hat{f}} }{L_k}}. 
\end{align}
Substituting \eqref{FGM_conv_eq_1} in the RHS of \eqref{t}, we obtain 
\begin{align}
	\label{kk}
	1 - \frac{\lambda_{k+1}}{\lambda_k} = \sqrt{ \frac{\mu_{\hat{f}}}{L_k} +\lambda_{k+1}\frac{\gamma_{0} - \mu_{\hat{f}}}{L_k}}. 
\end{align}
Making some algebraic manipulations to \eqref{kk} results in
\begin{align}
	\label{kkk}
	\frac{\lambda_k  - \lambda_{k+1}}{\lambda_k} &= \sqrt{\lambda_{k+1}} \sqrt{\frac{\mu_{\hat{f}}}{\lambda_{k+1} L_k} + \frac{\gamma_{0} - \mu_{\hat{f}}}{L_k}}.
\end{align}
Then, dividing both sides of \eqref{kkk} by $\lambda_{k+1}$, yields
\begin{align}
	\label{kkkkk}
	\frac{\lambda_k - \lambda_{k+1}}{\lambda_k \lambda_{k+1}} &= \frac{1}{\sqrt{\lambda_{k+1}}} \sqrt{\frac{\mu_{\hat{f}}}{\lambda_{k+1} L_k} + \frac{\gamma_{0} - \mu_{\hat{f}}}{L_k}}.
\end{align}
Note that we can re-write the LHS of \eqref{kkkkk} as
\begin{align}
	\frac{1}{\lambda_{k+1}} - \frac{1}{\lambda_k} &= \frac{1}{\sqrt{\lambda_{k+1}}}\sqrt{\frac{\mu_{\hat{f}}}{\lambda_{k+1} L_k} + \frac{\gamma_{0}  - \mu_{\hat{f}}}{L_k}}. 
\end{align}
Then, based on the difference of squares argument, we can derive that
	\begin{align}
		\label{convergence_stupid}
		\left(\frac{1}{\sqrt{\lambda_{k+1}}} - \frac{1}{\sqrt{\lambda_{k}}} \right) \left(\frac{1}{\sqrt{\lambda_{k+1}}} + \frac{1}{\sqrt{\lambda_{k}}}\right) &= \frac{1}{\sqrt{\lambda_{k+1}}} \sqrt{\frac{\mu_{\hat{f}}}{\lambda_{k+1} L_k} + \frac{\gamma_{0} - \mu_{\hat{f}}}{L_k}}.
	\end{align}

Let us now analyze the behavior of the terms in the sequence $\{ \lambda_k\}_{k=0}^\infty$. First, recall that from Lemma \ref{SFGM_lemma_2} we have $\alpha_k \in [0,1]$. Then, considering \eqref{lambda_recursive}, we can conclude that the terms $\lambda_k$ are non-increasing in the iteration counter $k$. Therefore, we can substitute the term $\frac{1}{\sqrt{\lambda_{k}}}$ in the LHS of \eqref{convergence_stupid} with the larger number $\frac{1}{\sqrt{\lambda_{k+1}}}$. This results in
\begin{align}
	\label{FGM_conv_eq_2}
	\frac{2}{\sqrt{\lambda_{k+1}}} \! \left( \! \frac{1}{\sqrt{\lambda_{k+1}}} - \frac{1}{\sqrt{\lambda_{k}}} \! \right)  &\geq \frac{1}{\sqrt{\lambda_{k+1}}} \sqrt{\frac{\mu_{\hat{f}}}{\lambda_{k+1} L_k} \! + \! \frac{\gamma_{0} \! - \! \mu_{\hat{f}}}{L_k}}
\end{align}

At this point, we note that the practical performance of the proposed method depends on the initialization of the parameter $\gamma_0$. To allow for the widest possible range of selection for this parameter, we need to consider separately the regions $\mathcal{R}_1 = [0, \mu_{\hat{f}}[$ and $\mathcal{R}_2 = [\mu_{\hat{f}}, 3L_k + \mu_{\hat{f}}]$. The results for the case when $\gamma_0 \in \mathcal{R}_2$ can be established by following the analysis conducted for FGM in \cite[Lemma 2.2.4]{Nesterov_book}. Therefore, in the sequel we will thoroughly analyze only the case when $\gamma_0 \in \mathcal{R}_1$, which is the novel part of the proposed proof. Let us begin by defining the following quantity 
\begin{align}
	\label{xi_k_def}
	\xi_{k, \mathcal{R}_1} \triangleq \sqrt{\frac{L_{\max}}{\left(\mu_{\hat{f}} - \gamma_{0} \right) \lambda_{k}}},
\end{align}
where $L_{\text{max}}$ was defined in \eqref{L_bound}. Next, we can rewrite \eqref{FGM_conv_eq_2} as
\begin{align}
	\label{sfgm_conv_useless}
	\frac{2}{\sqrt{\lambda_{k+1}}} - \frac{2}{\sqrt{\lambda_{k}}}  &\geq \sqrt{\frac{\mu_{\hat{f}} - \gamma_{0}}{L_k}} \sqrt{\frac{\mu_{\hat{f}} L_k}{L_k \lambda_{k+1} \left(\mu_{\hat{f}} - \gamma_{0} \right )} - 1}.
\end{align}
Then, we can relax the bound in \eqref{sfgm_conv_useless} and multiply it with a factor of $\sqrt{\frac{L_{\max}}{\mu_{\hat{f}} - \gamma_0}}$. This results in
\begin{align}
	\label{FGM_conv_eq_3}
	\xi_{k+1, \mathcal{R}_1} - \xi_{k, \mathcal{R}_1} &\geq \frac{1}{2}\sqrt{\frac{\mu_{\hat{f}} \xi_{k+1, \mathcal{R}_1}^2}{L_{\max}} - 1}.
\end{align}
We then proceed to establish via an inductive argument the following lower bound
\begin{align}
	\label{FGM_conv_eq_4}
	\xi_{k, \mathcal{R}_1} \geq \frac{1.959}{4 \delta} \sqrt{\frac{L_k}{\mu_{\hat{f}} - \gamma_0}} \left[e^{(k+1) \delta} - e^{(k+1) \delta}\right],
\end{align}
where $\delta \triangleq \frac{1}{2} \sqrt{\frac{\mu}{L}}$. Utilizing \eqref{xi_k_def} at step $k = 0$, we have
\begin{align}
	\label{loose_bound}
	\xi_{0, \mathcal{R}_1} &= \sqrt{\frac{L_{\max}}{(\mu_{\hat{f}} - \gamma_{0}) \lambda_{0}}} = \sqrt{\frac{L_{\max}}{\mu_{\hat{f}} - \gamma_{0}}},
\end{align}
where the second equality is obtained because $\lambda_{0} = 1$. Then, substituting \eqref{L_bound} in \eqref{loose_bound} we obtain 
\begin{align}
	\xi_{0, \mathcal{R}_1} &\geq \frac{1.959}{2} \sqrt{ \frac{L_k}{\mu_{\hat{f}} - \gamma_0}} \left[e^{1/2} - e^{-1/2}\right] \geq \frac{1.959}{4 \delta} \sqrt{ \frac{L_k}{\mu_{\hat{f}} - \gamma_0}} \left[e^\delta - e^{-\delta}\right] \label{tt}.
\end{align}
Note that the second row into \eqref{tt} follows because the RHS is increasing in $\delta$, which by construction is always $\delta < 0.5$.

As it is common with induction-type of proofs, the next step is to assume that \eqref{FGM_conv_eq_4} is satisfied for some iteration $k$. To establish that the relation would hold true at the next iteration as well, we proceed via contradiction. Define 
\[
\omega(t) \triangleq \frac{1.959}{4 \delta} \sqrt{\frac{L_k}{\mu_{\hat{f}} - \gamma_0}} \left[e^{(t+1) \delta} - e^{-(t+1) \delta}\right],
\] 
and note that from \cite[Lemma 2.2.4]{Nesterov_book} it is a convex function. Therefore, we can write
\begin{align}
	\label{FGM_conv_eq_5}
	\omega(t) \leq \xi_{k, \mathcal{R}_1} \stackrel{\eqref{FGM_conv_eq_3}}{\leq} \xi_{k+1, \mathcal{R}_1} - \frac{1}{2}\sqrt{\frac{\mu_{\hat{f}} \xi_{k+1, \mathcal{R}_1}^2}{L_{\max}} - 1}.
\end{align}
Assume that $\xi_{k+1, \mathcal{R}_1} < \omega(t+1)$. Then, substituting the definitions of $\delta$ and $\xi_{k, \mathcal{R}_1}$ into \eqref{FGM_conv_eq_5}, after massaging the resulting expression, we have
\begin{align}
	\omega(t) &< \omega(t+1) - \frac{1}{2}\sqrt{\frac{\mu_{\hat{f}} \xi_{k+1, \mathcal{R}_1}^2}{L_{\max}} - 1}.
\end{align}
Then, utilizing the definition of $\delta$, as well as \eqref{FGM_conv_eq_4}, we obtain
\begin{align}
	\omega(t) \! &\leq \! \omega(t \! + \! 1) \! - \frac{1}{2}\sqrt{4 \delta^2 \left[ \frac{1.959}{4 \delta} \sqrt{ \frac{L_k}{\mu_{\hat{f}} - \gamma_0}} \left(e^{(t+2) \delta} - e^{-(t+2) \delta}\right) \right]^2 \! \! \! \! - \! \! 1} \nonumber \\
	&= \omega(t+1) \!-\! \frac{1.959}{4} \sqrt{\!\frac{L_k}{\mu_{\hat{f}} - \gamma_0}} \left[e^{(t+2) \delta} \!+\! e^{-(t+2) \delta}\right] \!=\! \omega(t+1) \!+\! \omega'(t+1)\left(t \!-\! (t+1)\right) \leq \omega(t), \nonumber
\end{align}
where the last inequality follows from the supporting hyperplane theorem of convex functions. Notice that this result contradicts the earlier assumption that $\xi_{k+1, \mathcal{R}_1} < \omega(t+1)$. Thus, the inductive argument asserts that we have established the lower bound \eqref{FGM_conv_eq_4} to be true for all values of $k = 0, 1, \ldots$. 

We are finally ready to establish \eqref{FGM_conv_eq_66}. From \eqref{xi_k_def}, we can write
\begin{align}
	\label{jj}
	\lambda_{k} &= \frac{L_{\max}}{\xi_{k+1, \mathcal{R}_1}^2 (\mu_{\hat{f}} - \gamma_0)}. 
\end{align}
Utilizing \eqref{FGM_conv_eq_4} in the RHS of \eqref{jj}, we reach
\begin{align}
	\label{ttttttttttttt}
	\lambda_{k} \leq \frac{ (4 \delta)^2 L_{\max}}{1.959^2 L_k \left[e^{(k+1) \delta} - e^{(k+1) \delta}\right]^2},
\end{align}
All that remains in order to obtain the first inequality in \eqref{FGM_conv_eq_66} is to substitute definition of $\delta$ in \eqref{ttttttttttttt}.

To establish the remaining inequality in \eqref{FGM_conv_eq_66}, let us begin by analyzing the following quantity
\begin{align}
	\label{non-strongly-cvx-lambda}
	\left(e^{\frac{k + 1}{2} \sqrt{\frac{\mu_{\hat{f}}}{L_k}}} \! - \! e^{-\frac{k + 1}{2} \sqrt{\frac{\mu_{\hat{f}}}{L_k}}}\right)^2 \! &= \! e^{\left(k \! + \! 1\right) \sqrt{\frac{\mu_{\hat{f}}}{L_k}}} \! - \! e^{-\left(k+1\right) \sqrt{\frac{\mu_{\hat{f}}}{L_k}}} \! - \! 2. 
\end{align}
Then, utilizing the definition of the hyperbolic cosine function in \eqref{non-strongly-cvx-lambda}, we obtain
\begin{align}
	\left(e^{\frac{k + 1}{2} \sqrt{\frac{\mu_{\hat{f}}}{L_k}}} \! - \! e^{-\frac{k + 1}{2} \sqrt{\frac{\mu_{\hat{f}}}{L_k}}}\right)^2 \! = \! 2 \text{cosh}\left(\sqrt{\frac{\mu_{\hat{f}}}{L_k}} \left(k+1\right) \! - \! 2\right). 
\end{align}
From the Taylor expansion of the hyperbolic cosine function, we can write
\begin{align}
	\left(e^{\frac{k + 1}{2} \sqrt{\frac{\mu_{\hat{f}}}{L_k}}} - e^{-\frac{k + 1}{2} \sqrt{\frac{\mu_{\hat{f}}}{L_k}}}\right)^2 &= -2 + 2 + 2 \frac{\mu_{\hat{f}} \left(k+1\right)^2}{2L_k} + 2  \frac{\mu_{\hat{f}}^2 \left(k+1\right)^4}{4! {L_k}^2} + \ldots \text{.} \label{to be truncated}
\end{align}
The next step is to truncate the RHS of \eqref{to be truncated}. This results in
	\begin{align}
		\label{to be substituted}
		\left(e^{\frac{k + 1}{2} \sqrt{\frac{\mu_{\hat{f}}}{L_k}}} - e^{-\frac{k + 1}{2} \sqrt{\frac{\mu_{\hat{f}}}{L_k}}}\right)^2 \geq \frac{\mu_{\hat{f}}}{L_k} \left(k+1\right)^2 .
	\end{align}
All that remains to establish the second inequality of \eqref{FGM_conv_eq_66}, is to substitute \eqref{to be substituted} into the denominator of the first inequality of \eqref{FGM_conv_eq_66}. This completes the proof.
\section{Proof of Lemma 5}
\label{G}
We begin by substituting the upper bound \eqref{upper_bound} evaluated at the point $y = x^*$ into \eqref{FFF}, and obtain that
\begin{align}
	\label{F_0}
	F(x_0) &= \hat{f} (x_0) + \tau \hat{g} (x_0)
	\leq \hat{f} (x^*) \! + \! \nabla \hat{f} (x^*)^T \! (x_0 \! - \! x^*) \! + \! \frac{L_0}{2} ||x_0 \! - \! x^*||^2 \! + \! \tau \hat{g} (x_0).
\end{align}
Then, from the equality established in \eqref{reduced_grad}, we can write the RHS of \eqref{F_0} as
\begin{align}
	F(x_0) \leq &\hat{f} (x^*) \! + \! \nabla \hat{f} (x^*)^T \! (x_0 \! - \! x^*) \! + \! \frac{L_0}{2} ||x_0 \! - \! x^*||^2 \! + \! \tau \hat{g} (x_0) \nonumber \\ &= \hat{f} (x^*) \! + \!\left(\tau s_{L_0} (x^*) - L_0 \left(x^* - T_{L_0}(x^*)\right)\right)^T \! (x_0 \! - \! x^*) \!+ \frac{L_0}{2} ||x_0 - x^*||^2 + \tau \hat{g} (x_0). \label{22}
\end{align}
From the definition of the composite gradient mapping given in \eqref{T_L(y)}, we can see that when $y = x^*$, then $T_{L_0}(x^*) = x^*$. Therefore, we can now write the RHS of \eqref{22} as
\begin{align}
	F(x_0) &\leq \hat{f} (x^*) - \tau s_{L_0} (x^*) ^T (x^* - x_0) + \frac{L_0}{2} ||x_0 - x^*||^2 + \tau \hat{g} (x_0) . \label{222}
\end{align}
Lastly, utilizing \eqref{lower_bound_g} in the RHS of \eqref{222} completes the proof.

\end{document}

%% file: input.tex
\usepackage{amsfonts}
\usepackage{amsmath, amsthm}
\usepackage{amssymb}
\usepackage{cite}
\usepackage{dsfont}
\usepackage{latexsym}
\usepackage{times}
\usepackage{url}
\usepackage{verbatim}

\usepackage{graphicx}

\newtheorem{theorem}{Theorem}

\newtheorem{definition}{Definition}
\newtheorem{lemma}{Lemma}

\def\bb0{{\mathbb{0}}}

\def\bb{{\mathbf{b}}}

\def\b0{{\mathbf{0}}}

\def\sf0{{\mathsf{0}}}

\def\rm0{{\mathrm{0}}}


%% file: header.tex
\usepackage[]{hyperref}
\hypersetup{pdftex,colorlinks=true,allcolors=blue}
\usepackage{hypcap}
\usepackage{enumitem}
\usepackage[usenames, dvipsnames]{color}

